\documentclass[11pt]{amsart}
\usepackage{amsopn,amsmath,amsthm,amssymb}
\usepackage{latexsym}
\usepackage[hypertex]{hyperref}

\newcommand{\nc}{\newcommand}

\nc{\n}{\noindent}
\nc{\Id}{\mathbf{1}}
\nc{\vs}{\vspace{8pt}}
\nc{\alt}{\raise1pt\hbox{$\bigwedge$}}
\nc{\ncp}{\nabla^\mathrm{CP}}
\nc{\nhc}{\nabla^\mathrm{HC}}
\nc{\wncp}{\widehat\nabla^\mathrm{CP}}
\nc{\ct}{\cos_\theta}
\nc{\st}{\sin_\theta}
\nc{\ctt}{\cos_{\theta/2}}
\nc{\stt}{\sin_{\theta/2}}
\nc{\ft}{\hbox{$\frac12$}}
\nc{\Aff}{\mathit{Aff}}

\nc{\vg}{\mathfrak{v} }
\nc{\wg}{\mathfrak{w} }
\nc{\zg}{\mathfrak{z} }
\nc{\ngo}{\mathfrak{n} }
\nc{\kg}{\mathfrak{k} }
\nc{\mg}{\mathfrak{m} }
\nc{\bg}{\mathfrak{b} }
\nc{\ggo}{\mathfrak{g} }
\nc{\ggob}{\overline{\mathfrak{g}} }
\nc{\sog}{\mathfrak{so} }
\nc{\sug}{\mathfrak{su} }
\nc{\spg}{\mathfrak{sp} }
\nc{\slg}{\mathfrak{sl} }
\nc{\glg}{\mathfrak{gl} }
\nc{\cg}{\mathfrak{c} }
\nc{\hg}{\mathfrak{h} }
\nc{\tg}{\mathfrak{t} }
\nc{\ug}{\mathfrak{u} }
\nc{\dg}{\mathfrak{d} }
\nc{\ag}{\mathfrak{a} }
\nc{\pg}{\mathfrak{p} }
\nc{\sg}{\mathfrak{s} }
\nc{\pca}{\mathcal{P}}
\nc{\nca}{\mathcal{N}}

\nc{\ra}{\rightarrow}
\nc{\lra}{\longrightarrow}

\nc{\vp}{\varphi}
\nc{\ddt}{\frac{{\rm d}}{{\rm d}t}}

\nc{\SO}{{\tt SO}}
\nc{\Spe}{{\tt Sp}}
\nc{\SL}{{\tt SL}}
\nc{\SU}{{\tt SU}}
\nc{\Or}{{\tt O}}
\nc{\U}{{\tt U}}
\nc{\GL}{{\tt GL}}
\nc{\Se}{{\tt S}}
\nc{\CL}{{\tt CL}}
\nc{\Spin}{{\tt Spin}}
\nc{\Pin}{{\tt Pin}}

\nc{\RR}{{\mathbb R}}
\nc{\HH}{{\mathbb H}}
\nc{\CC}{{\mathbb C}}
\nc{\ZZ}{{\mathbb Z}}
\nc{\FF}{{\mathbb F}}
\nc{\NN}{{\mathbb N}}
\nc{\GG}{{\mathbb G}}
\nc{\JJ}{{\mathbb J}}
\nc{\II}{{\mathbb I}}
\nc{\KK}{{\mathbb K}}
\nc{\DD}{{\mathbb D}}
\nc{\EE}{{\mathbb E}}

\nc{\ad}{\operatorname{ad}}
\nc{\Ad}{\operatorname{Ad}}
\nc{\rank}{\operatorname{rank}}
\nc{\Irr}{\operatorname{Irr}}
\nc{\End}{\operatorname{End}}
\nc{\Aut}{\operatorname{Aut}}
\nc{\Inn}{\operatorname{Inn}}
\nc{\Der}{\operatorname{Der}}
\nc{\Ker}{\operatorname{Ker}}
\nc{\Iso}{\operatorname{I}}
\nc{\Le}{\operatorname{L}}
\nc{\tr}{\operatorname{tr}}
\nc{\dif}{\operatorname{d}\!}
\nc{\sen}{\operatorname{sen}}
\nc{\modu}{\operatorname{mod}}
\nc{\Ric}{\operatorname{R}}
\nc{\Sym}{\operatorname{Sym}}
\nc{\sca}{\operatorname{sc}}
\nc{\scalar}{{\sf s}}
\nc{\grad}{\operatorname{grad}}
\nc{\ric}{\operatorname{Ric}}
\nc{\Lie}{\operatorname{L}}
\nc{\tang}{\operatorname{T}}
\nc{\tm}{\operatorname{TM}}

\theoremstyle{plain}
\newtheorem{thm}{Theorem}[section]
\newtheorem{prop}[thm]{Proposition}
\newtheorem{cor}[thm]{Corollary}
\newtheorem{lem}[thm]{Lemma}

\theoremstyle{definition}

\theoremstyle{remark}
\newtheorem*{rem}{Remark}

\newtheorem{exam}[thm]{Example}

\newcommand{\ri}{{\rm (i)}}
\newcommand{\rii}{{\rm (ii)}}
\newcommand{\riii}{{\rm (iii)}}

\title[Complex product structures]{Complex product structures on 6-dimensional nilpotent Lie algebras}

\author{Adri\'an Andrada}
\address{The Abdus Salam International Centre for Theoretical Physics, Mathematics Section, Strada Cos\-tie\-ra 11, 34014 Trieste,
Italy}
\email{aandrada@ictp.trieste.it}
\address{CIEM, FaMAF, Universidad Nacional de C\'ordoba, Ciudad Universitaria, (5000) C\'or\-do\-ba, Argentina}
\email{andrada@mate.uncor.edu}

\begin{document}

\begin{abstract}
We study complex product structures on nilpotent Lie algebras,
establishing some of their main properties, and then we restrict
ourselves to 6 dimensions, obtaining the classification of
6-dimensional nilpotent Lie algebras admitting such structures. We prove that 
any complex structure which forms part of a complex product structure on a 6-dimensional nilpotent Lie algebra must be nilpotent in the sense of Cordero-Fern\'andez-Gray-Ugarte. 
A study is made of the torsion-free connection associated to the
complex product structure and we consider also the associated hypercomplex structures on the 12-dimensional nilpotent Lie algebras obtained by complexification.
\end{abstract}

\maketitle

\

\section{Introduction}
The notion of a complex product structure on a (real) Lie algebra
was introduced in \cite{AS}; it is given by two anticommuting
endomorphisms of the Lie algebra, one of them a complex structure
and the other one a product structure. Such a structure determines
in a unique way a torsion-free connection on the Lie algebra such
that the complex structure and the product structure are both
parallel, and this connection restricts to flat torsion-free
connections on two complementary totally geodesic subalgebras.

In \cite{AS} many examples of 4-dimensional Lie algebras admitting
complex product structures were exhibited, and the complete
classification of such Lie algebras was first given in \cite{BV}.
This classification, restricted to the solvable case, was also
shown in \cite{ABDO}. The next simplest step in order to obtain a
better understanding of these structures, would be the
classification of the 6-dimensional nilpotent Lie algebras which
admit them. The class of 6-dimensional
nilpotent Lie algebras has been extensively studied, and there is
only a finite number of such Lie algebras, up to isomorphism (see
for instance \cite{M}). Furthermore, the family of these Lie
algebras admitting a complex structure has also been specified,
and their list appears in \cite{S}.

Our goal in this work is to obtain the full classification of the
6-dimensional nilpotent Lie algebras admitting a complex product
structure, and this classification will be done independently of the
results given in \cite{S}. Our approach will make use of
properties of simply transitive actions of 3-dimensional nilpotent
Lie groups on 3-dimensional Euclidean space or, equivalently,
complete left symmetric algebra structures on 3-dimensional
nilpotent Lie algebras. The main result that will be used in the
classification is Theorem \ref{d-u}, which states that any
6-dimensional nilpotent Lie algebra admitting a complex product
structure can be decomposed as a direct sum (as vector spaces) of
a 4-dimensional subspace and a 2-dimensional central ideal, both
of them invariant by the complex and the product structures. As a
first corollary of this theorem we find two 6-dimensional
nilpotent Lie algebras admitting complex structures (according to
\cite{S}) that do not admit any complex product structure (see
Corollary \ref{non}). It will turn out later that the
6-dimensional nilpotent Lie algebras which admit complex
structures but admit no complex product structures are exactly 3.
Another corollary of Theorem \ref{d-u} is the fact that if a
complex structure on a 6-dimensional nilpotent Lie algebra forms
part of a complex product structure, then the complex structure
must be nilpotent in the sense of \cite{CFGU2} (see Corollary
\ref{Jnil}).

The outline of the paper is the following: \S 2 consists only of
necessary preliminaries, while in \S 3 we prove some general
results concerning complex product structures on nilpotent Lie
algebras. In \S 4 we restrict ourselves to the 6-dimensional case,
and we prove our main result (Theorem \ref{d-u}) which describes
the structure of such a Lie algebra, and this result will allow us
to perform, in \S 5, the classification mentioned earlier. In \S 6
we study the associated torsion-free connections, showing that
they are always complete and specifying whether these connections
are flat or not and, as an application, we use this information
together with results in \cite{AS} to obtain examples of
12-dimensional nilmanifolds equipped with hypercomplex structures.
The corresponding Obata connections are always Ricci flat, but not
necessarily flat. Finally, in \S 7 we consider some simple
examples at the Lie group level.

All Lie algebras and Lie groups in this work will be real and
finite dimensional, unless otherwise specified.

\

\section{Preliminaries}

Let us begin by recalling some basic facts. An {\em almost product
structure} on any Lie algebra $\ggo$ is a linear endomorphism
$E:\ggo \ra \ggo$ satisfying $E^2=\Id$ (and not equal to
$\pm\Id$). It is said to be {\em integrable} if
\[
E[x,y]=[Ex,y]+[x,Ey]-E[Ex,Ey] \quad \text{for
all } x,y\in\ggo.
\]
An integrable almost product structure will be called a {\em
product structure}. Given an almost product structure $E$ on
$\ggo$, we have a decomposition $\ggo=\ggo_+\oplus\ggo_-$ of
$\ggo$ into the direct sum of two linear subspaces, the
eigenspaces associated to the eigenvalues $\pm 1$ of $E$. It is
easy to verify that $E$ is integrable if and only if $\ggo_+$ and
$\ggo_-$ are both Lie subalgebras of $\ggo$; we will write
$\ggo=\ggo_+\bowtie\ggo_-$ and will call $(\ggo,\ggo_+,\ggo_-)$
the associated double Lie algebra. If $\dim \ggo_+=\dim\ggo_-$
then the product structure $E$ is called a {\em paracomplex
structure} (see \cite{KK, L}).

Next, we recall that an {\em almost complex structure} on any Lie
algebra $\ggo$ is a linear endomorphism $J:\ggo \ra \ggo$
satisfying $J^2=-\Id$. If $J$ satisfies the condition
\[
J[x,y]=[Jx,y]+[x,Jy]+J[Jx,Jy] \quad \text{for all } x,y\in \ggo,
\]
we will say that $J$ is {\em integrable} and we will call it a
{\em complex structure} on $\ggo$. Note that the dimension of a
Lie algebra carrying an almost complex structure must be even. A
complex structure on $\ggo$ induces a left-invariant complex
structure on $G$, the simply connected Lie group with Lie algebra
$\ggo$, so that $G$ becomes a complex manifold, although not
necessarily a complex Lie group.

If an almost complex structure $J$ satisfies the condition
$[Jx,Jy]=[x,y]$ for all $x,y\in\ggo$, then it is automatically
integrable, and it is called {\em abelian} \cite{BDM}. This name
follows from the fact that the eigenspaces corresponding to the
eigenvalues $\pm i$ of the natural extension
$J:\ggo^\CC\ra\ggo^\CC$ are abelian Lie subalgebras of the complex
Lie algebra $\ggo^\CC$. It has been proved in \cite{DF2} that only
solvable Lie algebras admit abelian complex structures.

Finally, a {\em complex product structure} on the Lie algebra
$\ggo$ is a pair $\{J,\,E\}$ of a complex structure $J$ and a
product structure $E$ satisfying $JE=-EJ$ (see \cite{AS}). If
$(\ggo,\ggo_+,\ggo_-)$ is the double Lie algebra associated to
$E$, then we have $\ggo_-=J\ggo_+$ and hence $E$ is in fact a
paracomplex structure. The endomorphism $F:=JE$ of $\ggo$ is
another product structure on $\ggo$ that anticommutes with $J$;
moreover, it can be seen that the operators in $\{\ct E+\st F:
\theta\in[0,2\pi)\}$ are all product structures anticommuting with
$J$. Two complex product structures $\{J,E\}$ and $\{J',E'\}$ on a
Lie algebra $\ggo$ are {\em equivalent} if there exists an
automorphism $\varphi$ of $\ggo$ such that $J'\varphi=\varphi J$
and $E'\varphi=\varphi E$.

If $\{J,E\}$ is a complex product structure on $\ggo$ and $J$ is
abelian, then it has been shown in \cite{AD} that the Lie
subalgebras $\ggo_+$ and $\ggo_-$ determined by $E$ are both
abelian, and the converse also holds. In this case, we will say
that $\{J,E\}$ is an {\em abelian} complex product structure.

Let us also recall that a complex product structure $\{J,E\}$ on
$\ggo$ determines uniquely a torsion-free connection $\ncp$ on
$\ggo$ such that $\ncp J=\ncp E=0$, i.e.
$\ncp:\ggo\times\ggo\ra\ggo$ is a $\ggo$-valued bilinear form on
$\ggo$ such that
\[ \ncp_xy-\ncp_yx=[x,y]\quad\text{(torsion-free)},\]
and also
\[ \ncp_xJy=J\ncp_xy,\quad \ncp_xEy=E\ncp_xy\quad\text{($J$ and $E$ are $\ncp$-parallel)}\]
for all $x,y\in\ggo$. From this we get that $\ncp_xy\in\ggo_{\pm}$
whenever $y\in\ggo_{\pm}$; in particular, $\ggo_+$ and $\ggo_-$
are totally geodesic subalgebras of $\ggo$ with respect to $\ncp$.
The induced connections $\nabla^+$ and $\nabla^-$ on $\ggo_+$ and
$\ggo_-$, respectively, are clearly torsion-free and also {\em
flat}, so that they define left symmetric algebra (LSA) structures
on these subalgebras. We recall that a {\em left-symmetric
algebra} structure on a Lie algebra $\hg$ is a bilinear product
$\hg\times\hg\longrightarrow\hg,\,(x,y)\mapsto x\cdot y$, which
satisfies the conditions
\begin{equation}\label{flat}
x\cdot(y\cdot z)-(x\cdot y)\cdot z=y\cdot(x\cdot z)-(y\cdot
x)\cdot z,
\end{equation}
\begin{equation}
[x,y]=x\cdot y-y\cdot x
\end{equation}
for all $x,y,z\in\hg$. The correspondence between a flat
torsion-free connection $\nabla$ and an LSA structure on $\hg$ is
given by $\nabla_xy=x\cdot y$, that is, $\nabla_x$ represents the
left-multiplication by $x\in\hg$. From \eqref{flat}, it follows
that $\nabla:\hg\rightarrow \glg(\hg)$ is a representation of
$\hg$. An LSA structure on $\hg$ is called {\em complete} if the
corresponding left-invariant connection on $H$, the simply
connected Lie group with Lie algebra $\hg$, is geodesically
complete. The completeness of an LSA structure can be
characterized in purely algebraic terms by the following
condition: the LSA structure is complete if and only if, for all
$x\in\hg$, the right multiplication $y\mapsto y\cdot x$ of $\hg$
is nilpotent (\cite{Se}).

We move on now to the Lie group level. Let $G$ be the simply
connected Lie group with Lie algebra $\ggo$, where we consider an
element of $\ggo$ as a left-invariant vector field on $G$. A
complex product structure $\{J,E\}$ on $\ggo$ determines a
left-invariant complex product structure on $G$, still denoted by
$\{J,E\}$. This means that $J$ is a complex structure and $E$ is a
product structure on the manifold $G$ such that $JE=-EJ$ and both
tensors are invariant by left-translations of $G$ (see \cite{A}
for a study of complex product structures on manifolds). If
$(\ggo,\ggo_+,\ggo_-)$ is the associated double Lie algebra, with
$\ggo_-=J\ggo_+$, let $G_+$ and $G_-$ denote the connected Lie
subgroups of $G$ with Lie algebras $\ggo_+$ and $\ggo_-$,
respectively. The decomposition $\ggo=\ggo_+\bowtie\ggo_-$
determines naturally two complementary distributions on $G$, both
of them involutive and the leaves of the foliations ${\mathcal
F}^+$ and ${\mathcal F}^-$ determined by these distributions are
totally real submanifolds of the complex manifold $(G,J)$.
Moreover, these leaves are totally geodesic and flat with respect
to the canonical torsion-free connection $\ncp$ determined by
$\{J,E\}$. It is easy to see that the leaf of ${\mathcal F}^{\pm}$
passing through $x\in G$ is $xG_{\pm}$, and it is well known that
this leaf is embedded if and only if $G_{\pm}$ is closed in $G$
(see for instance \cite{CL}).

\

\section{Complex product structures on nilpotent Lie algebras}

Given a Lie algebra $\ggo$, there is the associated lower central
series, which is defined recursively by
\[ \ggo^0=\ggo,\quad \ggo^k=[\ggo^{k-1},\ggo], \quad k\geq 1. \]
The Lie algebra $\ggo$ is called nilpotent if $\ggo^k=\{0\}$ for
some $k$. Equivalently, $\ggo$ is nilpotent if the endomorphisms
$\ad x:\ggo\ra\ggo$ are nilpotent for all $x\in\ggo$, where $(\ad
x) y=[x,y]$. If $\ggo^{k-1}\neq\{0\}$, but $\ggo^k=\{0\}$, then
$\ggo$ is said to be $k$-step nilpotent.

Recall that the Lie bracket on $\ggo$ can be thought of as a
linear map $[\,,]:\alt^2 \ggo \rightarrow \ggo$. In this way, we
can consider its transpose $\dif:\ggo^*\rightarrow \alt^2\ggo^*$,
defined as follows:
\[
(\dif f)(x\wedge y)= -f([x,y]) \quad \text{for all } f\in
\ggo^*,\, x,y\in \ggo.\] This linear mapping can be extended to
$\dif:\alt^k\ggo^*\ra\alt^{k+1}\ggo^*$ for $k\geq 1$, and the
Jacobi identity is equivalent to the vanishing of the composition
$\dif^2:\ggo^*\rightarrow \alt^3\ggo^*$. The nilpotency of a Lie
algebra $\ggo$ is equivalent to the existence of a basis
$\{e^1,\ldots,e^n\}$ of $\ggo^*$ such that
\[ \dif e^1=0,\quad \dif e^i\in\alt^2(\text{span}\{e^1,\ldots,e^{i-1}\}),\quad 2\leq i\leq n.\]
Following Salamon \cite{S}, we will use this characterization to
denote the isomorphism class of a nilpotent Lie algebra: we will
identify an $m$-dimensional Lie algebra with an $m$-tuple, where
in the $k^{\text{th}}$ spot we write $\dif e^k$, further
abbreviating $e^i\wedge e^j$ to $ij$. For instance, if we write
$\ggo=(0,0,0,0,12,14+23)$, this means that there is a basis
$\{e^1,\ldots,e^6\}$ of $\ggo^*$ such that $\dif e^i=0$ for $1\leq
i \leq 4$ and $\dif e^5=e^1\wedge e^2,\,\dif e^6=e^1\wedge
e^4+e^2\wedge e^3$. Equivalently, the dual basis
$\{e_1,\ldots,e_6\}$ of $\ggo$ satisfies
$[e_1,e_2]=-e_5,\,[e_1,e_4]=[e_2,e_3]=-e_6$.

General properties of nilpotent Lie algebras equipped with complex
structures have been given in \cite{S}. Moreover, in that paper the
full classification of 6-dimensional nilpotent Lie algebras admitting
complex structures was obtained. In forthcoming sections we will classify the
6-dimensional Lie algebras which carry a complex product
structure, but this classification will be made independently of the list
given by Salamon.

We mention next some properties of LSA structures on nilpotent Lie
algebras which will be needed later. We have already mentioned
that an LSA structure on any Lie algebra is complete if and only
if the right multiplication by any element of the algebra is
nilpotent. In the nilpotent case, we have the following equivalent
condition, proved by Kim:

\begin{thm}[\cite{K}]\label{kim}
An LSA structure on a nilpotent Lie algebra $\hg$ is complete if
and only if the left multiplication $y\mapsto x\cdot y$ is
nilpotent for all $x\in\hg$.
\end{thm}

Many nilpotent Lie algebras admit complete LSA structures; in
fact, it was even conjectured that every
solvable Lie algebra admits a complete LSA structure, until
Benoist \cite{B} gave the first example of a nilpotent Lie algebra
not admitting any. However, complete LSA structures on the
3-dimensional Heisenberg Lie algebra $\hg_3$ are well
understood, and we have the following result, proved by Fried and
Goldman:

\begin{thm}[\cite{FG}]\label{fg}
If $\nabla$ is a complete LSA structure on $\hg_3$, then
$\nabla_{z}=0$ for all $z\in\zg(\hg_3)$.
\end{thm}

It is well known that the existence of a complete LSA structure on
a Lie algebra $\ggo$ (which must be solvable according to a result
of Auslander \cite{Au}) is equivalent to a simply transitive
action of $G$, the simply connected Lie group with $Lie(G)=\ggo$,
on $\RR^{\dim G}$ as affine transformations. In this setting, the
theorem above implies that if the 3-dimensional Heisenberg group
acts simply transitively on $\RR^3$, then it contains a nontrivial
central translation.

\begin{rem}
Theorem \ref{fg} has been generalized in \cite{DDI} in the
following way: if a 2-step nilpotent Lie algebra $\ggo$ with
1-dimensional commutator ideal carries a complete LSA structure
$\nabla$, then there exists $z\in\zg(\ggo)$ such that
$\nabla_z=0$. It is known that a 2-step nilpotent Lie algebra with
1-dimensional commutator ideal is a central extension of a
Heisenberg Lie algebra $\hg_{2n+1}$ (for a proof see \cite{BD} or
\cite{DDI}). Recall that
$\hg_{2n+1}=\text{span}\{x_1,\ldots,x_n,y_1,\ldots,y_n,z\}$ with
Lie bracket given by $[x_j,y_j]=z$ for all $j=1,\ldots,n$.
\end{rem}

\medskip

We begin now the study of complex product structures on nilpotent
Lie algebras, focusing in the properties of the torsion-free
connection $\ncp$  associated to this structure. In Theorem
\ref{complete}, we prove that the LSA structures defined on the
subalgebras determined by the product structure are complete while
in Corollary \ref{ricci-flat} we show that $\ncp$ is Ricci-flat.

\begin{thm}\label{complete}
Let $\ggo$ be a nilpotent Lie algebra equipped with a complex
product structure $\{J,E\}$ and let $(\ggo,\ggo_+,\ggo_-)$ be the
associated double Lie algebra. Then, the flat torsion-free
connections $\nabla^+$ and $\nabla^-$ on $\ggo_+$ and $\ggo_-$,
respectively, determined by $\{J,E\}$ are complete.
\end{thm}

In order to prove this theorem, we will need the following general
result on product structures on nilpotent Lie algebras.

\begin{lem}\label{endo}
Let $E$ be a product structure on the nilpotent Lie algebra
$\ggo$, and let $(\ggo,\ggo_+,\ggo_-)$ be the associated double
Lie algebra. Define the representations
$\rho:\ggo_+\longrightarrow\glg(\ggo_-)$ and
$\mu:\ggo_-\longrightarrow\glg(\ggo_+)$ by
\begin{equation}\label{def}
[x,x']=-\mu(x')x+\rho(x)x',\end{equation}
for $x\in\ggo_+,\,x'\in\ggo_-$. Then $\rho(x)$ and $\mu(x')$ are
nilpotent endomorphisms of $\ggo_-$ and $\ggo_+$, respectively,
for all $x\in\ggo_+$ and $x'\in\ggo_-$.
\end{lem}

\begin{proof}
Let $\pi_+:\ggo\ra\ggo_+$ and
$\pi_-:\ggo\ra\ggo_-$ denote the projections. We will
prove inductively that
\begin{equation}\label{ro}
\pi_-\left((\ad x)^n(x')\right)=\rho(x)^n(x')
\end{equation}
for all $x\in\ggo_+,\,x'\in\ggo_-$. Clearly, \eqref{ro} holds for $n=1$. Next,
\begin{align*}
(\ad x)^{n+1}(x') & =[x,\pi_+\left((\ad x)^n(x')\right)+\rho(x)^n(x')]\\
                  & =[x,\pi_+\left((\ad x)^n(x')\right)]-\mu\left(\rho(x)^n(x')\right)x+\rho(x)^{n+1}(x').
\end{align*}
Since the first two terms in the last expression are in $\ggo_+$,
we get that the ($\ggo_-$)-component of $(\ad x)^{n+1}(x')$ is
$\rho(x)^{n+1}(x')$, and thus \eqref{ro} is proved.

In the same way one can prove that
\[
\pi_+\left((\ad x')^n(x)\right)=\mu(x')^n(x) \] for all
$x\in\ggo_+,\,x'\in\ggo_-$. Since the endomorphisms $(\ad x)$ and
$(\ad x')$ of $\ggo$ are nilpotent, then we obtain that
$\rho(x)\in\glg(\ggo_-)$ and $\mu(x')\in\glg(\ggo_+)$ are
nilpotent.
\end{proof}

\medskip

\begin{proof}[Proof of Theorem \ref{complete}]
We recall that $\nabla^+_{x_+}=-J\rho(x_+)J$ and
$\nabla^-_{x_-}=-J\mu(x_-)J$ for $x_+\in\ggo_+$ and
$x_-\in\ggo_-$, with $\rho$ and $\mu$ defined as in \eqref{def}
(see \cite{AS}). From Lemma \ref{endo}, we obtain at once that
$\nabla^+_{x_+}$ and $\nabla^-_{x_-}$ are nilpotent and thus both
$\nabla^+$ and $\nabla^-$ are complete, due to Theorem \ref{kim}.
\end{proof}

\medskip

We are going to show that the connection $\ncp$ associated to a
complex product structure on a nilpotent Lie algebra $\ggo$ is always
Ricci-flat. We recall that the Ricci tensor is defined by
$\ric(x,y)=\tr\left(z\mapsto R(z,x)y\right),\;x,y,z\in\ggo$, where
$R$ is the curvature tensor of $\ncp$.

\begin{prop}\label{tr}
Let $\{J,E\}$ be a complex product structure on the nilpotent Lie
algebra $\ggo$, and let $\ncp$ denote the associated torsion-free
connection. Then the endomorphisms $\ncp_z,\,z\in\ggo,$ are all
traceless.
\end{prop}

\begin{proof}
As $\ncp$ is parallel with respect to both $J$ and $E$, we know
that each endomorphism $\ncp_z$ of $\ggo$ must be of the form
$\begin{pmatrix} X & 0 \cr 0 & X \cr\end{pmatrix}$, where
$X\in\glg(n,\RR)$, with respect to suitable bases of $\ggo_+$ and
$\ggo_-$. Here $n$ is half the dimension of $\ggo$. If
$z=z_++z_-$, with $z_{\pm}\in\ggo_{\pm}$, then $X$ can be
considered as the matrix associated to the endomorphism
$(\nabla^+_{z_+}-J_-\nabla^-_{z_-}J_+)$ acting on $\ggo_+$,
where $J_+=J|_{\ggo_+}:\ggo_+\ra\ggo_-$ and
$J_-=J|_{\ggo_-}:\ggo_-\ra\ggo_+$. As $\nabla^+_{z_+}$ and
$\nabla^-_{z_-}$ are both nilpotent, they are traceless, and hence
so is $\ncp_z$ for all $z\in\ggo$.
\end{proof}

\begin{cor}\label{ricci-flat}
With the same hypothesis as in Proposition \ref{tr}, the Ricci
tensor $\ric$ of $\ncp$ vanishes.
\end{cor}

\begin{proof}
It has been proved in \cite{A} that $\ric$ is skew-symmetric and
that the following relation holds for any $x,y\in\ggo$:
\[ \ric(x,y)=-\ft\tr R(x,y).\]
Recalling that $R(x,y)=[\ncp_x,\ncp_y]-\ncp_{[x,y]}$, we have that
\[ \ric(x,y)=\ft\tr\ncp_{[x,y]}.\] From Proposition \ref{tr}, we
obtain that $\ric=0$.
\end{proof}

Let us move now to the Lie group level. So, let $G$ be the simply
connected nilpotent Lie group with Lie algebra $\ggo$, with
$\{J,E\}$ its left-invariant complex product structure. In this
case, it is known that $G$ is diffeomorphic to $\RR^{2n}$
($2n=\dim G$) via the exponential map $\exp:\ggo\lra G$ and that
the connected subgroups $G_+$ and $G_-$ corresponding to $\ggo_+$
and $\ggo_-$, respectively, are closed and simply connected, so that
the leaves of the complementary foliations ${\mathcal F}^{\pm}$
determined by them (the left cosets $xG_{\pm}$ for $x\in G$) are
all simply connected embedded submanifolds of $G$, diffeomorphic
to $\RR^n$. We already know that these leaves are totally geodesic
with respect to $\ncp$ and flat with respect to the induced
connection; furthermore, according to Theorem \ref{complete},
these leaves are also geodesically complete (this fact does not happen
in general).

A nilpotent Lie group $G$ admits a lattice (i.e., a cocompact
discrete subgroup) $\Gamma$ if and only if there is a basis of its
Lie algebra such that the corresponding structure constants are
rational (\cite{Ma}). If $G$ also admits a complex product
structure $\{J,E\}$, then it induces a complex product structure
$\{\tilde{J},\tilde{E}\}$ on the nilmanifold $M_\Gamma:=\Gamma\backslash G$,
and this complex product structure determines naturally two
complementary foliations on $M_\Gamma$, denoted by $\tilde{{\mathcal F}}^{\pm}$.
However, since $\Gamma$ is divided out from the left, the foliations
${\mathcal F}^{\pm}$ on $G$ determined by $\{J,E\}$ induce
foliations on $M_\Gamma$, denoted by ${\mathcal F}_{\Gamma}^{\pm}$.
We will show next that these foliations coincide.

\begin{prop}
Let $G$ be a Lie group equipped with a left-invariant complex
product structure $\{J,E\}$, and let ${\mathcal F}^{\pm}$ be the
associated complementary foliations on $G$. Suppose $\Gamma$ is a
lattice in $G$, and let $M_\Gamma:=\Gamma\backslash G$. Let
$\{\tilde{J},\tilde{E}\}$ be the complex product structure on
$M_\Gamma$ induced by $\{J,E\}$ and ${\mathcal F}_{\Gamma}^{\pm}$
be the foliations on $M_\Gamma$ induced by ${\mathcal F}^{\pm}$. If
$\tilde{{\mathcal F}}^{\pm}$ denote the foliations on $M_\Gamma$
determined by $\{\tilde{J},\tilde{E}\}$, then $\tilde{{\mathcal
F}}^{\pm}={\mathcal F}_{\Gamma}^{\pm}$.
\end{prop}

\begin{proof}
This proposition follows from the definition of the induced
objects in the quotient $M_\Gamma$, using the fact that the
projection $\pi:G\lra M_\Gamma$ is a covering map.
\end{proof}

\

\section{The 6-dimensional case}
We begin now the study of 6-dimensional nilpotent Lie algebras
admitting a complex product structure. In the main result of this
section we show that such a Lie algebra can be decomposed as the
direct sum (as vector spaces) of two subspaces invariant by the
complex product structure, one of which is a 2-dimensional ideal
contained in the centre of the algebra. This fact will allow us,
in the following section, to obtain the classification of
these algebras.

\medskip

\begin{thm}\label{d-u}
Let $\ggo$ be a 6-dimensional nilpotent Lie algebra equipped with
a complex product structure $\{J,E\}$. Then there is a
decomposition $\ggo=\dg\oplus\ug$ (direct sum of vector spaces),
where $\dg$ is a 4-dimensional subspace and $\ug$ is a
2-dimensional ideal contained in the centre $\zg$ of $\ggo$, such
that both $\dg$ and $\ug$ are invariant by $J$ and $E$.
\end{thm}

\begin{proof}
We know that $\ggo$ can be written as $\ggo=\ggo_+\bowtie\ggo_-$,
where $\ggo_-=J\ggo_+$. Since $\ggo_+$ and $\ggo_-$ are
3-dimensional nilpotent Lie algebras, they are isomorphic to
either the abelian Lie algebra $\RR^3$ or the Heisenberg Lie
algebra $\hg_3$. Recall that this latter Lie algebra has a basis
$\{e_1,e_2,e_3\}$ with Lie bracket given by $[e_1,e_2]=e_3$ and
$e_3$ is a central element. We will divide the study into two
cases: in the first one both subalgebras are abelian, whereas in
the second at least one of the subalgebras is isomorphic to
$\hg_3$.

\smallskip

{\bf (A)} First case: both $\ggo_+$ and $\ggo_-$ are abelian,
i.e., $\{J,E\}$ is an abelian complex product structure. We will
need the following result.

\begin{lem}\label{ab}
Let $\hg$ be a Lie algebra with an abelian complex product
structure and associated double Lie algebra
$(\hg,\hg_+,\hg_-),\,\hg_{\pm}\cong\RR^n$. If $\hg$ has a non
trivial centre $\zg$, then $\zg\cap\hg_+\neq\{0\}$,
$\zg\cap\hg_-\neq\{0\}$ and
$\zg=(\zg\cap\hg_+)\oplus(\zg\cap\hg_-)$. Moreover,
$J\left(\zg\cap\hg_+ \right)=\zg\cap\hg_-$.
\end{lem}

\begin{proof}
Let us suppose $\zg\cap\hg_+=\{0\}$, and take $0\neq z\in\zg$. We
can express $z=z_++z_-$ with $z_{\pm}\in\hg_{\pm},\,z_-\neq 0$.
Thus, $[z_+,x]=-[z_-,x]$ for all $x\in\hg$. Taking $x\in\hg_-$, we
have that $[z_+,x]=0$, and hence $z_+\in\zg\cap\hg_+$, so that
$z_+=0$. Therefore, $\zg\subset\hg_-$. However, since $J$ is
abelian, it preserves the centre of $\hg$, and thus $J\zg\subset
\zg\cap\hg_+=\{0\}$, which is impossible. Thus,
$\zg\cap\hg_+\neq\{0\}$, and the same holds analogously for $\hg_-$.
From the computations above it is clear that the equality
$\zg=(\zg\cap\hg_+)\oplus(\zg\cap\hg_-)$ holds.

Take now $z_+\in\zg\cap\hg_+$; since $J$ is an abelian complex
structure we have that, for $x_+\in\ggo_+$,
\[ [Jz_+,x_+]=-[z_+,Jx_+]=0.\]
This relation, together with the fact that $\hg_-$ is abelian,
implies that $Jz_+\in\zg\cap\hg_-$ and therefore
$J\left(\zg\cap\hg_+ \right)\subset\zg\cap\hg_-$. Reversing the
roles of $\hg_+$ and $\hg_-$ we obtain the opposite inclusion and
hence the proof of the lemma is complete.
\end{proof}

Applying Lemma \ref{ab}, we can find a basis $\{e_1,e_2,e_3\}$ of $\ggo_+$
and a basis $\{f_1,f_2,f_3\}$ of $\ggo_-$ such that
$f_i=Je_i,\,i=1,2,3$, and $e_3,f_3\in\zg$. Thus, the subspaces
$\dg=\text{span}\{e_1,e_2,f_1,f_2\}$ and
$\ug=\text{span}\{e_3,f_3\}$ satisfy the statement of the theorem.

\medskip

{\bf (B)} Second case: one of the subalgebras, say $\ggo_+$, is
isomorphic to $\hg_3$. We will prove first that the centre of
$\ggo_+$ is contained in the centre of $\ggo$. Let
$\{e_1,e_2,e_3\}$ be a basis of $\ggo_+$ such that $[e_1,e_2]=e_3$
and $[e_3,\ggo_+]=0$.

\begin{lem}\label{e-central}
The vector $e_3\in\zg(\ggo_+)$ belongs to the centre of $\ggo$.
\end{lem}

\begin{proof}
Let $x\in\ggo_-$. Then $[x,e_3]=\ncp_xe_3-\ncp_{e_3}x$, where
$\ncp$ is the torsion-free connection associated to the complex
product structure. Now,
\[ \ncp_{e_3}x=-J\ncp_{e_3}Jx=-J\nabla^+_{e_3}Jx=0 \]
due to Theorems \ref{fg} and \ref{complete}. Thus,
$[x,e_3]=\ncp_xe_3\in\ggo_+$. So we can write
\[ [x,e_3]=a_1e_1+a_2e_2+a_3e_3,\quad a_1,a_2,a_3\in\RR.\]
Now, on one hand $[[x,e_3],e_2]=a_1e_3$, but on the other hand
\[ [[x,e_3],e_2]=-[[e_3,e_2],x]-[[e_2,x],e_3]=-[[e_2,x],e_3]\]
using Jacobi's identity and $e_3\in\zg(\ggo_+)$. Hence,
$[[e_2,x],e_3]=-a_1e_3$, and since $\ggo$ is nilpotent, we have
$a_1=0$. Repeating this argument with $[[x,e_3],e_1]$, we arrive
at $a_2=0$ and so we get $[x,e_3]=a_3e_3$, which implies that
$a_3=0$. Therefore $[x,e_3]=0$ and $e_3\in\zg(\ggo)$.
\end{proof}

\medskip

We will consider now two subcases, according to the isomorphism
class of $\ggo_-$.

\smallskip\n $\ri$ Let us suppose first that $\ggo_-$ is abelian.
If we define $f_i\in\ggo_-$ by $f_i=Je_i,\,i=1,2,3$, then
$\{f_1,f_2,f_3\}$ is a basis of $\ggo_-$. Using Lemma
\ref{e-central}, we arrive at the following result.

\begin{lem}\label{f3}
The element $f_3\in\ggo_-$ belongs to the centre of $\ggo$.
\end{lem}

\begin{proof}
Clearly, $[f_3,\ggo_-]=0$. So, let us take now $x\in\ggo_+$ and
perform the following computations, using the integrability of
$J$:
\[ J[f_3,x]=-[e_3,x]+[f_3,Jx]-J[e_3,Jx]=0\]
since $e_3$ is central and $\ggo_-$ is abelian. Therefore, $f_3$
is also central.
\end{proof}
The subspaces $\dg=\text{span}\{e_1,e_2,f_1,f_2\}$ and
$\ug=\text{span}\{e_3,f_3\}$ satisfy the statement of the theorem.

\smallskip\n $\rii$ Assume now that $\ggo_-$ is isomorphic to
$\hg_3$, and let us consider a basis $\{f_1,f_2,f_3\}$ of $\ggo_-$
such that $[f_1,f_2]=f_3$ and $[f_3,\ggo_-]=0$. From Lemma
\ref{e-central} applied to both $\ggo_+$ and $\ggo_-$, we know
that both $e_3$ and $f_3$ belong to the centre of $\ggo$ and from
this we get the following result.

\begin{lem}\label{h3x2}
With notation as above, we may suppose that $Je_1=f_1,\,Je_2=f_2$
and $Je_3=cf_3$ for some $c\neq 0$.
\end{lem}

\begin{proof}
Let the complex structure $J$ be given by
\[ \begin{cases}
Je_1=\alpha_1f_1+\alpha_2f_2+\alpha_3f_3,\\
Je_2=\beta_1f_1+\beta_2f_2+\beta_3f_3,\\
Je_3=\gamma_1f_1+\gamma_2f_2+\gamma_3f_3,
\end{cases}\]
with $\alpha_i,\beta_i,\gamma_i\in\RR,\,i=1,2,3,$ and let us
suppose that $\gamma_1^2+\gamma_2^2\neq 0$. By interchanging $f_1$
with $f_2$, if necessary, we can assume further that $\gamma_1\neq
0$. By subtracting from $e_1$ and $e_2$ a suitable multiple of
$e_3$, we obtain a new `$e_1$' and a new `$e_2$' such that
$[e_1,e_2]=e_3$ and
\[ \begin{cases}
Je_1=\alpha_2f_2+\alpha_3f_3,\\
Je_2=\beta_2f_2+\beta_3f_3,\\
Je_3=\gamma_1f_1+\gamma_2f_2+\gamma_3f_3,
\end{cases} \]
with $\gamma_1\neq 0$ and
$\Delta:=\alpha_2\beta_3-\alpha_3\beta_2\neq 0$. Note that
$Jf_3=\Delta^{-1}(\beta_2e_1-\alpha_2e_2)$. From the integrability
of $J$, we have
\[ Je_3=J[e_1,e_2]=[Je_1,e_2]+[e_1,Je_2]+J[Je_1,Je_2]=-\alpha_2[e_2,f_2]+\beta_2[e_1,f_2]\]
so that
\[ Je_3=\beta_2[e_1,f_2]-\alpha_2[e_2,f_2].\]
Using again the integrability of $J$, we get
\[ 0=J[e_1,e_3]=[e_1,Je_3]+J[Je_1,Je_3]=\gamma_1[e_1,f_1]+\gamma_2[e_1,f_2]-\alpha_2\gamma_1Jf_3\]
and
\[ 0=J[e_2,e_3]=[e_2,Je_3]+J[Je_2,Je_3]=\gamma_1[e_2,f_1]+\gamma_2[e_2,f_2]-\beta_2\gamma_1Jf_3\]
so that
\[ \begin{cases}
\gamma_1[e_1,f_1]+\gamma_2[e_1,f_2]=\alpha_2\gamma_1Jf_3,\\
\gamma_1[e_2,f_1]+\gamma_2[e_2,f_2]=\beta_2\gamma_1Jf_3.
\end{cases} \]
Now, let us compute
\begin{align*}
[[e_1,f_2],Je_3] & =[[e_1,f_2],\gamma_1f_1+\gamma_2f_2+\gamma_3f_3] \\
                 & =\gamma_1[[e_1,f_1],f_2]+\gamma_2[[e_1,f_2],f_2]\\
                 & =[\gamma_1[e_1,f_1]+\gamma_2[e_1,f_2],f_2]\\
                 & =\alpha_2\gamma_1[Jf_3,f_2]\\
                 & =\frac{\alpha_2\gamma_1}{\Delta}\left(\beta_2[e_1,f_2]-\alpha_2[e_2,f_2]\right)\\
                 & =\frac{\alpha_2\gamma_1}{\Delta}Je_3,
\end{align*}
and since $\ggo$ is nilpotent and $\gamma_1\neq0$, we must have
$\alpha_2=0$. Let us compute next
\begin{align*}
[[e_2,f_2],Je_3] & =[[e_2,f_2],\gamma_1f_1+\gamma_2f_2+\gamma_3f_3] \\
                 & =\gamma_1[[e_2,f_1],f_2]+\gamma_2[[e_2,f_2],f_2]\\
                 & =[\gamma_1[e_2,f_1]+\gamma_2[e_2,f_2],f_2]\\
                 & =\beta_2\gamma_1[Jf_3,f_2]\\
                 & =\frac{\beta_2\gamma_1}{\Delta}Je_3,
\end{align*}
and again since $\ggo$ is nilpotent and $\gamma_1\neq 0$, we get
that $\beta_2=0$. Therefore, we have that both $Je_1$ and $Je_2$
are in the subspace spanned by $f_3$, which is impossible since
they must be linearly independent. Therefore, the condition
$\gamma_1^2+\gamma_2^2\neq 0$ can never hold, and as a consequence
$Je_3$ is a multiple of $f_3$. Since $[Je_1,Je_2]=c^{-1}f_3$, for
some $c\neq 0$, we can take $\{Je_1,Je_2,c^{-1}f_3\}$ as the new
basis of $\ggo_-$, which satisfies the conditions stated in the
lemma.
\end{proof}

Resuming the proof of the theorem, we see that the subspaces
$\dg=\text{span}\{e_1,e_2,f_1,f_2\}$ and
$\ug=\text{span}\{e_3,f_3\}$ satisfy the conditions in the
statement, and therefore the proof is now complete.
\end{proof}

\medskip

As a first consequence of this theorem, together with results in
\cite {S}, we have the following

\begin{cor}\label{non}
The Lie algebras $(0,0,0,12,23,14-35)$ and $(0,0,12,13,23,14+25)$
possess complex structures, but do not admit any complex product
structure.
\end{cor}

\begin{proof}
These Lie algebras carry complex structures, as shown in \cite{S}.
However, they cannot admit any complex product structure since
each centre is one dimensional, and thus there cannot exist a
2-dimensional ideal $\ug$ contained in the centre.
\end{proof}

\smallskip

Let us recall now the definition of a {\em nilpotent complex
structure} on a nilpotent Lie algebra, which was introduced in
\cite{CFGU1}. Given a complex structure $J$ on a nilpotent Lie
algebra $\ggo$, there is the ascending series $\{\ag_t(J)\}$
associated to $J$, which is defined inductively by
\[ \ag_0(J)=\{0\},\quad \ag_t(J)=\left\{x\in\ggo:
[x,\ggo]\subset\ag_{t-1}(J),\;[Jx,\ggo]\subset\ag_{t-1}(J)\right\},\]
for $t\geq 1$. $J$ is called {\em nilpotent} if $\ag_t(J)=\ggo$
for some $t\geq 1$. Compact nilmanifolds equipped with nilpotent
complex structures (i.e. induced by a nilpotent complex structure
on the corresponding Lie algebra) possess many interesting
properties (see \cite{CFGU2}). As another corollary of Theorem
\ref{d-u}, we have the following:

\begin{cor}\label{Jnil}
If $\{J,E\}$ is a complex product structure on a
6-dimensional nilpotent Lie algebra, then $J$ is a nilpotent
complex structure.
\end{cor}

\begin{proof}
In \cite{CFGU1} it was proved that
a complex structure $J$ on a 6-dimensional nilpotent Lie algebra
$\ggo$ is nilpotent if and only if $\ag_1(J)\neq \{0\}$; note that
$\ag_1(J)=\{x\in\zg(\ggo):Jx\in\zg(\ggo)\}$.

Consider now the complex product structure $\{J,E\}$
on $\ggo$. According to Theorem \ref{d-u}, there exists a central
ideal $\ug$ in $\ggo$ which is $J$-invariant. Thus,
$\ag_1(J)\supset\ug\neq \{0\}$ and hence $J$ is nilpotent.
\end{proof}

\medskip

In general, the dimension of the centre of a nilpotent Lie algebra admitting
a complex product structure need not be $\geq 2$ and, furthermore,
the complex structure need not be nilpotent, as in the
6-dimensional case. The next example will illustrate these facts.

\begin{exam}
We exhibit next an example of an 8-dimensional nilpotent Lie
algebra with 1-dimensional centre which admits a complex product
structure. Consider the 4-dimensional nilpotent Lie algebra
$\ngo_4$ with a basis $\{e_1,\ldots,e_4\}$ such that
$[e_1,e_2]=e_3,\,[e_1,e_3]=e_4$, and the complete LSA structure $\nabla$
on $\ngo_4$ given by:
\[\nabla_{e_1}=\begin{pmatrix} 0&0&0&0 \cr 0&0&0&0 \cr 0&1&0&0 \cr
0&0&1&0 \cr\end{pmatrix},\quad \nabla_{e_2}=\begin{pmatrix}
0&0&0&-1 \cr 0&0&0&0 \cr 0&0&0&0 \cr 0&0&0&0 \cr\end{pmatrix},\]
\[  \nabla_{e_3}=\begin{pmatrix} 0&0&1&0 \cr 0&0&0&0 \cr 0&0&0&0 \cr
0&0&0&0 \cr \end{pmatrix}, \quad \nabla_{e_4}=\begin{pmatrix}
0&-1&0&0 \cr 0&0&0&0 \cr 0&0&0&0 \cr 0&0&0&0 \cr \end{pmatrix},\]
in the ordered basis $\{e_1,\ldots,e_4\}$. Observe that
$e_4\in\zg(\ngo_4)$ but $\nabla_{e_4}\neq 0$ (compare Theorem
\ref{fg}); this LSA structure on $\ngo_4$ was given by Fried in
\cite{F} (see also \cite{K}). Let $V$ denote the 4-dimensional
vector space underlying $\ngo_4$, so that
$\nabla:\ngo_4\rightarrow \glg(V)$ defines a representation of
$\ngo_4$ on $V$ and hence we can form the semidirect product
$\ggo:=\ngo_4\ltimes_{\nabla} V$. Due to results in \cite{AS}, the
Lie algebra $\ggo$ admits a complex product structure $\{J,E\}$
whose associated double Lie algebra is $(\ggo,\ngo_4,V)$. Note
that in this case the central element $e_4$ of $\ngo_4$ does not
belong to the centre of $\ggo$; furthermore,
$\zg(\ggo)=\text{span}\{(0,e_1)\}$ is 1-dimensional. It can also
be seen that the complex structure $J$ satisfies $\ag_t(J)=\{0\}$
for all $t\geq 0$; therefore, the complex structure $J$ is not
nilpotent.
\end{exam}

\smallskip

\begin{rem}
Not every nilpotent complex structure on a 6-dimensional nilpotent
Lie algebra may be part of a complex product structure. For
instance, consider the Lie algebra $\ggo=(0,0,0,0,13+42,14+23)$
and the nilpotent complex structure $J$ on $\ggo$ defined by
$Je_1=e_2,\,Je_3=e_4,\,Je_5=e_6$. This complex structure satisfies
$J[x,y]=[Jx,y]$ for all $x,y\in\ggo$, and it has been shown in
\cite{AS} that a complex structure satisfying such condition
cannot be part of a complex product structure, unless the algebra
is abelian. Note that this Lie algebra does admit complex product
structures (see \S5), where the complex structure is different
from the one considered here.
\end{rem}

\

\section{Classification}

Let $\{J,E\}$ be a complex product structure on the 6-dimensional
nilpotent Lie algebra $\ggo$. It is known that if $J$ is abelian,
then the associated double Lie algebra must be of the form
$(\ggo,\RR^3,\RR^3)$. On the other hand, if $J$ is not abelian,
then there are essentially two possibilities for the associated
double Lie algebra: either $(\ggo,\hg_3,\RR^3)$ or
$(\ggo,\hg_3,\hg_3)$. However, in the next result we will show
that a Lie algebra admits a complex product structure of the
former type if and only if it admits a complex product structure
of the latter type. This fact will simplify notoriously the task
of classifying the 6-dimensional nilpotent Lie algebras admitting
complex product structures, which will be carried out next.

\begin{prop}\label{je}
Let $\{J,E\}$ be a complex product structure on the 6-dimensional
nilpotent Lie algebra $\ggo$ with $(\ggo,\ggo_+,\ggo_-)$ the
associated double Lie algebra.

$\ri$ If $\ggo_+\cong\hg_3,\,\ggo_-\cong\hg_3$, then there exists
a complex product structure $\{J,E'\}$ on $\ggo$ such that its
associated double Lie algebra is $(\ggo,\hg_3,\RR^3)$, where
$E'=\ct E+\st JE$ for some $\theta\in [0,2\pi)$.

$\rii$ If $\ggo_+\cong\hg_3,\,\ggo_-\cong\RR^3$, then there exists
a complex product structure $\{J,E'\}$ on $\ggo$ such that its
associated double Lie algebra is $(\ggo,\hg_3,\hg_3)$, where $E'=JE$.
\end{prop}

\begin{proof}

$\ri$ According to Lemma \ref{h3x2}, there exists a basis
$\{e_1,e_2,e_3\}$ of $\ggo_+$ and a basis $\{f_1,f_2,f_3\}$ of
$\ggo_-$ such that $[e_1,e_2]=e_3,\,[f_1,f_2]=f_3$ (with some non
zero $[e_i,f_j]$), $Je_i=f_i$ for $i=1,2$ and $Je_3=cf_3$ with
$c\neq 0$. Consider the following subspaces of $\ggo$:
\begin{gather*}
\ggo'_+=\text{span}\left\{ce_1+f_1,\,ce_2+f_2,\,e_3+f_3\right\},\\
\ggo'_-=\text{span}\left\{e_1-cf_1,\,e_2-cf_2,\,-\frac{1}{c}e_3+cf_3\right\}.
\end{gather*}
The complementary subspaces $\ggo'_+$ and $\ggo'_-$ are in fact
Lie subalgebras of $\ggo$; moreover,
$\ggo'_+\cong\hg_3,\;\ggo'_-\cong\RR^3$ and $J\ggo'_+=\ggo'_-$.
The associated product structure $E'$ is given by $E'=\ct
E+\st JE$, where
\[ \cos\theta=\frac{c^2-1}{c^2+1},\quad
\sin\theta=\frac{2c}{c^2+1},\] and $\ri$ is proved.

$\rii$ There exists a basis $\{e_1,e_2,e_3\}$ of $\ggo_+$ and a
basis $\{f_1,f_2,f_3\}$ of $\ggo_-$ such that $[e_1,e_2]=e_3$
(with some non zero $[e_i,f_j]$) and $Je_i=f_i$ for $i=1,2,3$.
Consider the following subspaces of $\ggo$:
\begin{gather*}
\ggo'_+=\text{span}\left\{e_1+f_1,\,e_2+f_2,\,e_3+f_3\right\},\\
\ggo'_-=\text{span}\left\{e_1-f_1,\,e_2-f_2,\,e_3-f_3\right\}.
\end{gather*}
The complementary subspaces $\ggo'_+$ and $\ggo'_-$ are in fact
Lie subalgebras of $\ggo$; moreover,
$\ggo'_+\cong\hg_3,\;\ggo'_-\cong\hg_3$ and $J\ggo'_+=\ggo'_-$.
The associated product structure $E'$ is given by $E'=JE$ and
hence $\rii$ is proved.
\end{proof}

\medskip

Let us consider the Lie algebra $\tilde{\ggo}:=\ggo/\ug$, where we
use the notation from Theorem \ref{d-u}. It is a 4-dimensional
nilpotent Lie algebra, and since $\ug$ is invariant under $J$ and
$E$, it carries an induced complex product structure
$\{\tilde{J},\tilde{E}\}$. There are only two nilpotent Lie
algebras of dimension 4 which admit complex product structures,
the abelian one $\RR^4$ and the central extension of the
Heisenberg Lie algebra $\hg_3\times\RR$. We consider next each of
these cases separately.

\subsection{First case: $\tilde{\ggo}=\hg_3\times\RR$}\label{h3r}

In \cite{AS} there is a classification of complex product
structures on $\hg_3\times\RR$, which is given as follows. Let
$\{v_1,v_2,v_3,v_4\}$ be a basis of $\hg_3\times\RR$ such that
$[v_1,v_2]=v_3$ and $v_3,v_4$ are central elements. Then every
complex product structure on this Lie algebra is equivalent to
$\{J',E_\theta\}$, where the complex structure $J'$ is
$J'v_1=v_2,\,J'v_3=v_4$ and the product structure $E_{\theta}$ is
given by
\[ E_{\theta}=\begin{pmatrix} 1 & 0 & 0 & 0 \cr 0 & -1 & 0 & 0 \cr 0 & 0 & \cos\theta &
                      \sin\theta \cr 0 & 0 & \sin\theta & -\cos\theta \cr \end{pmatrix} \]
for some $\theta\in[0,2\pi)$, in the ordered basis
$\{v_1,v_2,v_3,v_4\}$. Any automorphism of $\tilde{\ggo}$ which
defines an equivalence between $\{\tilde{J},\tilde{E}\}$ and
$\{J',E_{\theta}\}$ can be lifted to an automorphism of $\ggo$
which is the identity on $\ug$, and via this latter automorphism
we can construct a new complex product structure on $\ggo$
equivalent to the original one, and such that the induced complex
structure on $\tilde{\ggo}$ is $\{J',E_{\theta}\}$. Therefore, we
can suppose without loss of generality that $\tilde{J}=J'$ and
$\tilde{E}=E_{\theta}$ for some $\theta\in [0,2\pi)$. The
subalgebras of $\hg_3\times\RR$ associated to the eigenvalues $\pm
1$ of $\tilde E$ are
\[ \tilde{\ggo}_+=\text{span}\{v_1,\,\ctt v_3+\stt v_4\},\]
and \[ \tilde{\ggo}_-=\text{span}\{v_2,\,-\stt v_3+\ctt v_4\}.\]
If $p:\ggo\lra\tilde{\ggo}$ denotes the canonical projection, then
we have $p(\ggo_+)=\tilde{\ggo}_+$ and $p(\ggo_-)=\tilde{\ggo}_-$;
furthermore, $p(\dg\cap\ggo_+)=\tilde{\ggo}_+$ and
$p(\dg\cap\ggo_-)=\tilde{\ggo}_-$ as vector spaces. Hence, we can
find a basis $\{e_1,e_2,e_3\}$ of $\ggo_+$ and a basis
$\{f_1,f_2,f_3\}$ of $\ggo_-$ such that $\dg=\text{span}\{e_1,e_2,f_2,f_2\}$,
$\ug=\text{span}\{e_3,f_3\}$, $f_i=Je_i$ for $i=1,2,3$, and
\begin{equation}\label{original} \begin{cases}
[e_1,e_2]=\alpha e_3\,(\alpha\in\{0,1\}),\\
[f_1,f_2]=\beta f_3\,(\beta\in\{0,1\}),\\
[e_1,f_1]=Ae_2+Be_3+Cf_2+Df_3,\,(A^2+C^2\neq 0)\\
[e_1,f_2]\in\ug,\\
[e_2,f_1]\in\ug,\\
[e_2,f_2]\in\ug.
\end{cases}
\end{equation}
Also, from the integrability of $J$ we get that
\[ J[e_1,e_2]=[f_1,e_2]+[e_1,f_2]+J[f_1,f_2],\]
and so
\begin{equation}\label{alpha}
[e_1,f_2]-[e_2,f_1]=\beta e_3+\alpha f_3.
\end{equation}
Let us denote $[e_2,f_2]=Ge_3+Hf_3$ for some $G,H\in\RR$. By Jacobi's
identity, $[[e_1,f_1],e_2]=-[[f_1,e_2],e_1]-[[e_2,e_1],f_1]=0$
since $\ug\subset\zg$. On the other hand,
$[[e_1,f_1],e_2]=-C[e_2,f_2]$, so that $CG=CH=0$. In the same way
we verify that $0=[[e_1,f_1],f_2]=A[e_2,f_2]$, and hence
$AG=AH=0$. As $A^2+C^2\neq 0$, we must have $G=H=0$, so that
\[[e_2,f_2]=0.\]

We have now two subcases: $\ri\;\alpha=\beta=0$ and
$\rii\;\alpha=1,\,\beta=0$. In fact, the case $\alpha=0,\,\beta=1$
is analogous to the case $\rii$ above, and the case
$\alpha=1,\,\beta=1$ need not be considered due to Proposition
\ref{je}.

\medskip

\subsubsection{$\alpha=0,\,\beta=0$.}
In this case the relations \eqref{original}, under condition
\eqref{alpha}, become simply
\begin{equation}\label{100} \begin{cases} [e_1,f_1]=Ae_2+Be_3+Cf_2+Df_3,\\
[e_2,f_1]=[e_1,f_2]=Ee_3+Ff_3.\end{cases} \end{equation} If
$E=F=0$, then $\ggo\cong(0,0,0,0,0,12)$. Let us consider now
$E^2+F^2\neq 0$. Suppose first $A\neq 0$. Defining
\begin{gather*}
v_1:=Ae_1+Cf_1,\;v_2:=-f_1,\;v_3:=Fe_3-Ef_3,\\
v_4:=Af_2,\;v_5:=A[e_1,f_1],\;v_6:=-A^2[e_2,f_1],\end{gather*} we
have that $\{v_1,\ldots,v_6\}$ is a basis for $\ggo$ and
\[ [v_1,v_2]=-v_5,\quad [v_2,v_5]=-v_6,\quad [v_1,v_4]=-v_6.\]
So, $\ggo\cong(0,0,0,0,12,14+25)$. If $A=0$, so that $C\neq 0$, we
define
\[ v_1:=Cf_1,\;v_2:=e_1,\;v_3:=Fe_3-Ef_3 ,\;
v_4:=-Ce_2,\;v_5:=C[e_1,f_1],\;v_6:=-C[e_1,f_2],  \] and with this
basis we see again that $\ggo\cong(0,0,0,0,12,14+25)$. This
concludes this case.

\medskip

\subsubsection{$\alpha=1,\,\beta=0$.} The relations
\eqref{original}, under condition \eqref{alpha}, become in this
case \begin{equation}\label{110}
\begin{cases}
[e_1,e_2]=e_3,\\
[e_1,f_1]=Ae_2+Be_3+Cf_2+Df_3,\,(A^2+C^2\neq 0)\\
[e_2,f_1]=Ee_3+Ff_3,\\
[e_1,f_2]=Ee_3+(F+1)f_3.
\end{cases} \end{equation}
Changing $f_1$ by $f_1-Be_2$ we may assume that $B=0$, so that
$[e_1,f_1]=Ae_2+Cf_2+Df_3$.

Let us suppose first $E\neq 0$ and $AF-CE\neq0$. In this case, the
vectors $v_1,\ldots,v_6$ defined below form a basis of $\ggo$:
\begin{gather*} v_1:=E^{-1}f_1,\;v_2:=Ee_1+(F+1)f_1,\;v_3:=Ee_2+Ff_2,\\
v_4:=[e_1,f_1],\;v_5:=[e_2,f_1],\;v_6:=(AF-CE)[e_1,f_2]. \end{gather*}
These vectors satisfy the relations
\[ [v_1,v_2]=-v_4,\, [v_1,v_3]=-v_5,\,[v_1,v_4]=-AE^{-1}v_5,\,[v_2,v_4]=-v_6.  \]
If $A=0$, we see immediately that $\ggo\cong(0,0,0,12,13,24)$. If
$A\neq 0$, multiplying $v_3$ and $v_5$ by $AE^{-1}$, we see that
$\ggo\cong(0,0,0,12,13+14,24)$. Let us still suppose that
$E\neq0$, but now with $AF-CE=0$; note that $A\neq 0$. Consider
the vectors
\begin{gather*} v_1:=-(Ee_1+(F+1)f_1),\; v_2:=E^{-1}f_1,\; v_3:=E^{-1}f_2,\\
v_4:=[e_1,f_1],\;v_5:=[e_1,f_2],\;v_6:=AE^{-1}[e_2,f_1].\end{gather*}
Then we have
\[ [v_1,v_2]=-v_4,\;[v_1,v_3]=-v_5,\;[v_2,v_4]=-v_6,\]
so that $\ggo\cong(0,0,0,12,13,24)$.

Let us now move to the case $E=0$. Consider first $F=0$, and
define
\[ v_1:=-e_1,\;v_2:=f_1,\;v_3:=e_2,\;v_4:=[e_1,f_1],\;v_5:=e_3,\;v_6:=Ae_3+Cf_3.\]
These vectors form a basis of $\ggo$ if $C\neq 0$, and as they
satisfy
\[ [v_1,v_2]=-v_4,\;[v_1,v_3]=-v_5,\;[v_1,v_4]=-v_6,\]
we have that $\ggo\cong(0,0,0,12,13,14)$. If $C=0$ (and hence
$A\neq 0)$, we consider the basis
\[ v_1:=-e_1,\;v_2:=f_1,\;v_3:=f_2,\;v_4:=[e_1,f_1],\;v_5:=f_3,\;v_6:=Ae_3,\]
and we have again that $\ggo\cong(0,0,0,12,13,14)$. Consider now
$F=-1$, so that
\[ [e_1,e_2]=e_3,\quad [e_1,f_1]=Ae_2+Cf_2+Df_3,\quad [e_2,f_1]=-f_3.\]
If $A=0$, then $[e_1,f_1]\in\zg$ and it is easily seen that
$\ggo\cong(0,0,0,12,13,23)$. If $A\neq 0$, define
\[ v_1:=-e_1,\;v_2:=f_1,\;v_3:=f_2,\;v_4:=[e_1,f_1],\;v_5:=Ae_3,\;v_6:=-Af_3.\]
Then, $[v_1,v_2]=-v_4,\;[v_1,v_4]=-v_5,\;[v_2,v_4]=-v_6$ and thus,
$\ggo\cong (0,0,0,12,14,24)$. Let us finally suppose that $F\neq
0,\,F\neq -1$. If $A\neq 0$, consider the basis
\[ v_1:=-e_1,\;v_2:=f_1,\;v_3:=-\frac{AF}{F+1}f_2,\]
\[ v_4:=[e_1,f_1],\;v_5:=Ae_3+C(F+1)f_3,\;v_6:=-AFf_3.\]
We have then
\[ [v_1,v_2]:=-v_4,\; [v_1,v_4]=-v_5,\;[v_1,v_3]=-v_6,\;[v_2,v_4]=v_6,\]
and thus $\ggo\cong (0,0,0,12,14,13+42)$. If $A=0$ (and hence
$C\neq 0$), we consider the basis
\[ v_1:=-e_1,\;v_2:=f_1,\;v_3:=\frac{C(F+1)}{F}e_2,\]
\[v_4:=[e_1,f_1],\;v_5:=\frac{C(F+1)}{F}e_3,\;v_6:=C(F+1)f_3, \]
and we have
\[ [v_1,v_2]=-v_4,\; [v_1,v_3]=-v_5,\;[v_1,v_4]=-v_6,\;[v_2,v_3]=-v_6,\]
so that $\ggo\cong(0,0,0,12,13,14+23)$. This concludes this case.

\medskip

\subsection{Second case: $\tilde{\ggo}=\RR^4$}\label{r4}
Any complex product structure on $\RR^4$ is equivalent to
$\{J_0,E_0\}$, where $J_0$ and $E_0$ are given by
\[ J_0=\begin{pmatrix} 0 & -\Id\cr \Id & 0
\cr \end{pmatrix},\quad E_0=\begin{pmatrix} \Id & 0\cr 0 & -\Id\cr
\end{pmatrix}\]
in some ordered basis of $\RR^4$, where $\Id$ is the $(2\times
2)$-identity matrix. If $\{\tilde{J},\tilde{E}\}$ is the complex
product structure on $\tilde{\ggo}$ induced by $\{J,E\}$, we can
suppose without loss of generality that $\tilde{J}=J_0$ and
$\tilde{E}=E_0$. If $p:\ggo\lra\tilde{\ggo}$ denotes the canonical
projection, then we have $p(\ggo_+)=\tilde{\ggo}_+$ and
$p(\ggo_-)=\tilde{\ggo}_-$; furthermore,
$p(\dg\cap\ggo_+)=\tilde{\ggo}_+$ and
$p(\dg\cap\ggo_-)=\tilde{\ggo}_-$ as vector spaces. Hence, we can
find a basis $\{e_1,e_2,e_3\}$ of $\ggo_+$ and a basis
$\{f_1,f_2,f_3\}$ of $\ggo_-$ such that $\dg=\text{span}\{e_1,e_2,f_1,f_2\}$,
$\ug=\text{span}\{e_3,f_3\}$, $f_i=Je_i$ for $i=1,2,3$, and
\begin{equation}\label{original2} \begin{cases}
[e_1,e_2]=\alpha e_3\,(\alpha\in\{0,1\}),\\
[f_1,f_2]=\beta f_3\,(\beta\in\{0,1\}),\\
[e_i,f_j]\in\ug,\quad 1\leq i,j\leq 2.\\
\end{cases}
\end{equation}
Also, from the integrability of $J$ we get that
\[ J[e_1,e_2]=[f_1,e_2]+[e_1,f_2]+J[f_1,f_2],\]
and so
\begin{equation}\label{beta}
[e_1,f_2]-[e_2,f_1]=\beta e_3+\alpha f_3.
\end{equation}
Note that, since $\tilde{\ggo}$ is abelian, the commutator ideal
$[\ggo,\ggo]\subseteq\ug$ and therefore $\ggo$ is 2-step nilpotent
with $\dim[\ggo,\ggo]\leq 2$.

We have now two subcases: $\ri\;\alpha=\beta=0$ and
$\rii\;\alpha=1,\,\beta=0$. In fact, the case $\alpha=0,\,\beta=1$
is analogous to the case $\rii$ above, and the case
$\alpha=1,\,\beta=1$ need not be considered due to Proposition
\ref{je}.

\medskip

\subsubsection{$\alpha=0,\,\beta=0$.} In this case, the relations
\eqref{original2}, under the condition \eqref{beta}, become \[
\begin{cases}
[e_1,f_1]=A_1e_3+A_2f_3,\\
[e_1,f_2]=[e_2,f_1]=B_1e_3+B_2f_3,\\
[e_2,f_2]=D_1e_3+D_2f_3.
\end{cases}\]
Performing computations as the ones done in \S\ref{h3r}, which we
omit, we arrive at the following result:

\n $\diamond$ If $\dim[\ggo,\ggo]=1$, then $\ggo$ is isomorphic
either to $(0,0,0,0,0,12)$ or $(0,0,0,0,0,12+34)$.

\n $\diamond$ If $\dim[\ggo,\ggo]=2$, then $\ggo$ is isomorphic to
one of the following Lie algebras: \[(0,0,0,0,12,34), \,
(0,0,0,0,13+42,14+23), \, (0,0,0,0,12,14+23).\]

\medskip

\subsubsection{$\alpha=1,\,\beta=0$.} In this case, the relations
\eqref{original2}, under the condition \eqref{beta}, become
\[ \begin{cases}
[e_1,e_2]=e_3,\\
[e_1,f_1]=A_1e_3+A_2f_3,\\
[e_2,f_1]=C_1e_3+C_2f_3,\\
[e_1,f_2]=C_1e_3+(C_2+1)f_3,\\
[e_2,f_2]=D_1e_3+D_2f_3.
\end{cases} \]
Note that in this case we have $\dim[\ggo,\ggo]=2$. Computing
again as in \S\ref{h3r}, we can conclude that $\ggo$ must be
isomorphic to one of the following Lie algebras:
\[(0,0,0,0,12,13),\,(0,0,0,0,13+42,14+23),\,(0,0,0,0,12,14+23), (0,0,0,0,12,34). \]

\medskip

\subsection{Conclusion}\label{conc}
We summarize in the following table the results obtained in the
previous subsections.

\

\begin{center}
\small{
\begin{tabular}{|l|c|c|c|}\hline
 & $\RR^3\bowtie\RR^3$ & $\hg_3\bowtie\RR^3$ & $\hg_3\bowtie\hg_3$\\ \hline

$(0,0,0,0,0,0)$ & yes & no & no \\ \hline

$(0,0,0,0,0,12)$ & yes & no  & no \\ \hline

$(0,0,0,0,0,12+34)$ & yes & no & no\\ \hline

$(0,0,0,0,12,14+25)$ & yes & no & no\\ \hline

$(0,0,0,0,12,13)$ & no & yes & yes \\ \hline

$(0,0,0,0,13+42,14+23)$ & yes & yes & yes \\ \hline

$(0,0,0,0,12,14+23)$ & yes & yes & yes \\ \hline

$(0,0,0,0,12,34)$ & yes & yes & yes \\ \hline

$(0,0,0,12,13,14)$ & no & yes & yes \\ \hline

$(0,0,0,12,13,23)$ & no & yes & yes \\ \hline

$(0,0,0,12,14,24)$ & no & yes & yes \\ \hline

$(0,0,0,12,13,24)$ & no & yes & yes \\ \hline

$(0,0,0,12,13+14,24)$ & no & yes & yes \\ \hline

$(0,0,0,12,13,14+23)$ & no & yes & yes \\ \hline

$(0,0,0,12,14,13+42)$ & no & yes & yes \\ \hline
\end{tabular} }
\end{center}

\

\begin{rem}
Let us identify some of the Lie algebras appearing in the list
above. The Lie algebra $(0,0,0,0,0,12)$ is the product
$\hg_3\times\RR^3$; the Lie algebra $(0,0,0,0,0,12+34)$ is the
product $\hg_5\times\RR$, while the Lie algebra $(0,0,0,0,12,34)$
is the product $\hg_3\times\hg_3$. The algebra
$(0,0,0,0,13+42,14+23)$ is isomorphic to the complex 3-dimensional
Heisenberg algebra $\hg_3^{\CC}$, considered as a real Lie
algebra. Finally, the Lie algebra $(0,0,0,12,13,23)$ is isomorphic
to the 2-step nilpotent free Lie algebra on 3 generators. Indeed,
this Lie algebra is isomorphic to the vector space
$\RR^3\oplus\alt^2\RR^3$, with Lie bracket given only by
$[u,v]=u\wedge v$ for $u,v\in\RR^3$.
\end{rem}

\medskip

We already know that there are two 6-dimensional nilpotent Lie
algebras which admit a complex structure but no complex product
structure, see Corollary \ref{non}. However, comparing with the
list in \cite{S}, we see that there is another algebra that admits
a complex structure but no complex product structure, as it does
not appear in the table above; this Lie algebra is the one whose
Lie bracket is encoded in $(0,0,0,12,13+42,14+23)$. Next, we
verify this fact by direct computations.

\begin{prop}
The Lie algebra $\ggo=(0,0,0,12,13+42,14+23)$ does not admit any
complex product structure.
\end{prop}

\begin{proof}
There is a basis $\{e_1,\ldots,e_6\}$ of $\ggo$ which satisfies:
\[ [e_1,e_2]=-e_4,\,[e_1,e_3]=-e_5,\,[e_2,e_4]=e_5,\,[e_1,e_4]=[e_2,e_3]=-e_6.\]
Note that $\zg(\ggo)=\text{span}\{e_5,e_6\}$.

Let us suppose that $\ggo$ admits a complex product structure
$\{J,E\}$ with associated double Lie algebra $(\ggo,\ggo_+,\ggo_-),\,
\ggo_-=J\ggo_+$. It is known that $\ggo$ admits both abelian and
non abelian complex structures (see \cite{CFU,S}). In the former
case we have that the associated double Lie algebra is isomorphic
to $(\ggo,\RR^3,\RR^3)$, whereas in the latter case it is
isomorphic to either $(\ggo,\hg_3,\RR^3)$ or $(\ggo,\hg_3,\hg_3)$.
However, according to Proposition \ref{je}, we may simply suppose
that it is isomorphic to $(\ggo,\hg_3,\RR^3)$. Therefore in any
case we have that at least one of the subalgebras $\ggo_{\pm}$ is
abelian; we may assume hence that $\ggo_-$ is abelian. From Lemmas
\ref{ab} and \ref{f3}, we can choose a basis $\{x,y,z\}$ of
$\ggo_-$ with $z\in\zg(\ggo)$. If we write $x=\sum x_ie_i,\,y=\sum
y_ie_i$ with $x_i,y_i\in\RR,\,i=1,\ldots,6$, then from $[x,y]=0$
we obtain the equations
\begin{equation}\label{g1} \begin{cases}
x_1y_2-x_2y_1=0,\\
x_1y_3-x_3y_1-x_2y_4+x_4y_2=0,\\
x_1y_4-x_4y_1+x_2y_3-x_3y_2=0.
\end{cases} \end{equation}

Suppose first $x_1=x_2=y_1=y_2=0$; thus $e_1,e_2\in\ggo_+$ and
since $\ggo_+$ is a subalgebra, we have that
$e_4=-[e_1,e_2]\in\ggo_+$. But then $e_6=-[e_1,e_4]$ is also in
$\ggo_+$, which contradicts the fact that $\dim\ggo_+=3$. Hence,
$x_1^2+x_2^2+y_1^2+y_2^2\neq 0$. Now, from the first equation in
\eqref{g1}, we have that $(y_1,y_2)=\alpha
(x_1,x_2),\,\alpha\in\RR$, where we can assume without loss of
generality that $x_1^2+x_2^2\neq 0$. Then, substituting $y$ by
$y-\alpha x$, we are allowed to consider $y_1=y_2=0$. Using this
in \eqref{g1}, we get
\[ \begin{cases}
x_1y_3-x_2y_4=0,\\
x_1y_4+x_2y_3=0.
\end{cases} \]
As $x_1^2+x_2^2\neq 0$, we have that $y_3=y_4=0$, and hence
$y\in\zg(\ggo)$. So, $\zg(\ggo)$ is contained in $\ggo_-$, but
this contradicts Theorem \ref{d-u}. Therefore, there is no complex
product structure on $\ggo$.
\end{proof}

\begin{rem}
The Lie algebra $\ggo=(0,0,0,12,13+42,14+23)$ is an example of a
Lie algebra that admits both complex and product structures, but
does not admit any complex product structure. In fact, the almost
complex structure $J$ on $\ggo$ given by
\[ Je_1=e_2,\quad Je_3=-e_4,\quad Je_5=e_6,\]
is integrable, (in fact, it is abelian, see \cite{CFU}). Also,
$\ggo$ admits paracomplex structures. One example of such a
structure is given by the following decomposition of $\ggo$ as the
sum of two 3-dimensional subalgebras $\ggo=\ggo_+\oplus\ggo_-$,
where
\[ \ggo_+=\text{span}\{e_1,e_3,e_5\},\quad
\ggo_-=\text{span}\{e_1+e_2,e_3+e_4,e_6\}.\]
\end{rem}

\

\section{Associated torsion-free connections}

In this section we will determine when the torsion-free connection
$\ncp$ associated to a complex product structure on the Lie
algebras classified before are flat, that is, they are LSA
structures on these Lie algebras. In what follows, we will use the
following notation: as $\ncp$ is parallel with respect to both $J$
and $E$, then each endomorphism $\ncp_x$ will be of the form
$\begin{pmatrix} X & 0 \cr 0 & X\cr \end {pmatrix}$, where
$X\in\glg(3,\RR)$, with respect to suitable bases of $\ggo_+$ and
$\ggo_-$, so that we will simply write $\ncp=X$. We recall also
the definition of $\ncp$, which is given as follows:
\begin{equation}\label{casos}
  \begin{cases}
    \ncp_{x_+}y_+= -\pi_+J[x_+,Jy_+],\\
    \ncp_{x_-}y_-= -\pi_-J[x_-,Jy_-],\\
    \ncp_{x_+}y_-= \pi_-[x_+,y_-],\\
    \ncp_{x_-}y_+= \pi_+[x_-,y_+],
  \end{cases}
\end{equation}
for any $x_+,y_+\in\ggo_+,\;x_-,y_-\in\ggo_-$, where
$\pi_{\pm}:\ggo\ra\ggo_{\pm}$ are the projections (see
\cite{A,AS}).

First, let us state a general result on the torsion-free
connection associated to a complex product structure $\{J,E\}$,
when the product structure $E$ is replaced by another product
structure in $\{\ct E+\st F: \theta\in[0,2\pi)\}$. This result
holds also for complex product structures on manifolds, with
exactly the same proof.

\begin{prop}
Let $\{J,E\}$ be a complex product structure on $\ggo$ and let
$\nabla:=\ncp$ be its associated torsion-free connection. Let
$E_\theta:=\ct E+\st JE$; if $\nabla^\theta$ denotes the
torsion-free connection associated to the complex product
structure $\{J,E^\theta\}$, then $\nabla^\theta=\nabla$.
\end{prop}

\begin{proof}
Let us recall that $\nabla^\theta$ is the only torsion-free
connection on $\ggo$ such that $\nabla^\theta J=\nabla^\theta
E_\theta=0$. Let us compute now the following:
\begin{multline*}
\nabla_x E_\theta y=\nabla_x(\ct E+\st JE)y=\ct\nabla_x Ey+\st\nabla_x JEy=\\
=\ct E\nabla_x y+\st JE\nabla_x y=(\ct E+\st JE)\nabla_xy=E_\theta
\nabla_xy, \end{multline*} so that $\nabla E_\theta=0$. As also
$\nabla J=0$, we obtain that both connections coincide.
\end{proof}

As a consequence, from this Proposition together with Proposition
\ref{je} we obtain that we have to consider only the complex
product structures whose associated double Lie algebra is either
$(\ggo,\RR^3,\RR^3)$ or $(\ggo,\hg_3,\RR^3)$. We will do this in
each of the two main cases: when $\tilde{\ggo}=\hg_3\times\RR$ or
when $\tilde{\ggo}=\RR^4$. Let us begin now studying each case
separately.

\subsection{First case: $\tilde{\ggo}=\hg_3\times\RR$}\

\n $\ri$ $\alpha=0,\;\beta=0$. In this case we can see that
\[ \ncp_{e_1}=\begin{pmatrix} 0&0&0 \cr C&0&0 \cr D&F&0 \cr \end{pmatrix},\quad
\ncp_{e_2}=\begin{pmatrix} 0&0&0 \cr 0&0&0 \cr F&0&0 \cr \end{pmatrix}, \]
\[ \ncp_{f_1}=\begin{pmatrix} 0&0&0 \cr -A&0&0 \cr -B&-E&0 \cr \end{pmatrix},\quad
\ncp_{f_2}=\begin{pmatrix} 0&0&0 \cr 0&0&0 \cr -E&0&0 \cr
\end{pmatrix} \] and $\ncp_{e_3}=\ncp_{f_3}=0$. Therefore, it is
readily verified that
\[ R(e_1,f_1)e_1=-2(AF-CE)e_3,\;R(e_1,f_1)f_1=-2(AF-CE)f_3,\]
and all the other possibilities equal to zero. Thus, $\ncp$ is
flat if and only if
\[ AF=CE.\]

\smallskip

\n $\rii$ $\alpha=1,\;\beta=0$. In this case, we see that
\[ \ncp_{e_1}=\begin{pmatrix} 0&0&0 \cr C&0&0 \cr D&F+1&0 \cr \end{pmatrix},\quad
\ncp_{e_2}=\begin{pmatrix} 0&0&0 \cr 0&0&0 \cr F&0&0 \cr \end{pmatrix}, \]
\[ \ncp_{f_1}=\begin{pmatrix} 0&0&0 \cr -A&0&0 \cr -B&-E&0 \cr \end{pmatrix},\quad
\ncp_{f_2}=\begin{pmatrix} 0&0&0 \cr 0&0&0 \cr -E&0&0
\cr \end{pmatrix} \]
and $\ncp_{e_3}=\ncp_{f_3}=0$. One can easily see that
\[ R(e_1,f_1)e_1=-(2(AF-CE)+A)e_3,\;R(e_1,f_1)f_1=-(2(AF-CE)+A)f_3,\]
so that $\ncp$ is flat if and only if
\[ A(2F+1)=2CE.\]

\medskip

\subsection{Second case: $\tilde{\ggo}=\RR^4$}\
In this case we can prove that $\ncp$ is always flat, irrespective
of the values of $\alpha$ and $\beta$.

\begin{prop}\label{rflat}
If $\tilde{\ggo}\cong\RR^4$, then the torsion-free connection
$\ncp$ on $\ggo$ associated to the complex product structure on
$\ggo$ is flat.
\end{prop}

\begin{proof}
From the equations \eqref{original2} and the definition of $\ncp$
(see \eqref{casos}), we can easily deduce that
\[ \ncp_xy\in\ug \quad\text{for all } x,y\in\ggo.\]
Also, it is easy to verify that
\[ \ncp_ux=\ncp_xu=0 \quad\text{for all } u\in\ug,\,x\in\ggo.\]
These facts, together with $[\ggo,\ggo]\subset\ug$, imply that
$\ncp$ is flat.
\end{proof}

\smallskip

\begin{rem}
In the case $\tilde{\ggo}=\RR^4,\;\alpha=\beta=0$, we could have
proved the flatness of $\ncp$ using the following general result:

\begin{prop}
Any abelian complex product structure on a 2-step nilpotent Lie
algebra is flat, i.e. the associated torsion-free connection
$\ncp$ is flat.
\end{prop}

\begin{proof}
Although this result might be proved in a straightforward manner,
we will give an indirect proof. From results in \cite{AS}, we know
that any complex product structure on a Lie algebra $\ggo$ gives
rise to a hypercomplex structure on $(\ggo^\CC)_\RR$. Furthermore,
if the complex product structure is abelian, then the hypercomplex
structure is also abelian (meaning that each complex structure is
abelian). It is also known that $\ncp$ (the torsion-free
connection on $\ggo$ associated to the complex product structure)
is flat if and only if $\nhc$ (the Obata connection on
$(\ggo^\CC)_\RR$ associated to the hypercomplex structure) is
flat. Now, the proposition follows since in \cite{DF} it has been
proved that the Obata connection associated to an abelian
hypercomplex structure on a 2-step nilpotent Lie algebra is flat.
\end{proof}
\end{rem}

\

Combining the results obtained so far in this section with the
results from \S\ref{h3r} and \S\ref{r4}, we obtain

\begin{thm} \label{nonflat}
$\ri$ The following Lie algebras admit only flat complex product
structures:

$\begin{array}{lcl}
(0,0,0,0,0,0), & \quad & (0,0,0,0,12,14+23), \\
(0,0,0,0,0,12), & \quad & (0,0,0,0,12,34), \\
(0,0,0,0,0,12+34), & \quad & (0,0,0,12,13,23), \\
(0,0,0,0,12,13), & \quad & (0,0,0,12,13,14+23). \\
(0,0,0,0,13+42,14+23), & \quad & \\
\end{array}$

\smallskip

\n $\rii$ The following Lie algebras admit only non flat complex
product structures:

$\begin{array}{l}
(0,0,0,12,14,24), \\ (0,0,0,12,13,24).
\end{array}$

\smallskip

\n $\riii$ The following Lie algebras admit both flat and non flat
complex product structures:

$\begin{array}{lcl}
(0,0,0,0,12,14+25), & \quad & (0,0,0,12,13+14,24), \\ (0,0,0,12,13,14), & \quad & (0,0,0,12,14,13+42).
\end{array}$
\end{thm}

\medskip

We shall consider now the question of the completeness of these
connections and we will prove that they are always complete. We
recall first that a left-invariant connection on a Lie group $H$
is complete if and only if the following differential equation on
its Lie algebra $\hg$ admits solutions defined on the whole real
line: \begin{equation}\label{completa}
\dot{x}(t)=-\nabla_{x(t)}x(t),\end{equation} where $\nabla$
denotes the bilinear form on $\hg$ induced by the connection (see
for instance \cite{BM}).

\begin{prop}
If $\ncp$ is the torsion-free connection associated to a complex
product structure on a 6-dimensional nilpotent Lie algebra, then
$\ncp$ is complete.
\end{prop}

\begin{proof}
We will consider again different cases, according to the ones in
\S\ref{h3r} and \S\ref{r4}.

\smallskip {\bf First case: $\tilde{\ggo}=\hg_3\times\RR$.}

\n$\ri$ $\alpha=0,\,\beta=0$. Suppose that the curve $x(t)$ in
equation \eqref{completa} is given by $x(t)=\sum a_i(t)e_i+\sum
b_i(t)f_i$, with $a_i,\,b_i$ real valued functions defined on some
interval of the real line. From the brackets given in \eqref{100}
and the expressions for $\ncp$ given earlier in this section, we
obtain that the differential equation \eqref{completa} yields the
system
\[ \begin{cases}
\dot{a}_1=0,\\
\dot{a}_2=-Ca_1^2+Aa_1b_1,\\
\dot{a}_3=-Da_1^2-2Fa_1a_2+Ba_1b_1+Ea_2b_1+Ea_1b_2,\\
\dot{b}_1=0,\\
\dot{b}_2=-Ca_1b_1+Ab_1^2,\\
\dot{b}_3=-Da_1b_1-Fa_1b_2-Fa_2b_1+Bb_1^2-2Eb_1b_2.
\end{cases} \]
From this it follows immediately that $a_1$ and $b_1$ are constant
functions, $a_2,b_2$ are linear functions and $a_3,b_3$ are
quadratic functions, all of them defined on the whole real line.
Thus, the connection $\ncp$ in this case is complete.

\n$\rii$ $\alpha=1,\,\beta=0$. From the brackets given in
\eqref{110} and the expressions for $\ncp$ given early in this
section, we obtain that the differential equation \eqref{completa}
yields a system similar to the one above, and it is easy to verify
that its solutions are all polynomials of degree $\leq 2$, and
hence the connection is complete. We omit the details.

\medskip {\bf Second case: $\tilde{\ggo}=\RR^4$.} We know from
Proposition \ref{rflat} that $\ncp$ is flat, so that it defines an
LSA structure on $\ggo$. This LSA structure is complete if and
only the right multiplications $\rho(x)$ are nilpotent or,
equivalently, if $\tr\rho(x)=0$ for all $x\in\ggo$ (see
\cite{Se}). Since $\ncp$ is torsion-free, we have that
$\ad(x)=\ncp_x-\rho(x)$ for all $x\in\ggo$, but $\ad(x)$ is
nilpotent and $\ncp_x$ is traceless, so that $\tr\rho(x)=0$ for
all $x\in\ggo$, and therefore $\ncp$ is complete.
\end{proof}

\medskip

\subsection{Application: 12-dimensional hypercomplex nilpotent Lie
algebras}

In \cite{AS} it was proved that if a Lie algebra $\ggo$ carries a
complex product structure, then its complexification considered as
a real Lie algebra, i.e. $(\ggo^\CC)_\RR$, is endowed with a
hypercomplex structure. It was shown also that if $\ncp$ stands
for the torsion-free connection on $\ggo$ associated to the
complex product structure, its natural extension to
$(\ggo^\CC)_\RR$ coincides with $\nhc$, the Obata connection
associated to this hypercomplex structure. Therefore, $\nhc$ is
flat if and only $\ncp$ is flat (see \cite{AS}).

Applying this to the 6-dimensional nilpotent Lie algebras equipped
with a complex product structure, we obtain a list of
12-dimensional hypercomplex nilpotent Lie algebras. The associated
Obata connections may be flat or non flat, depending on the
flatness of the connection associated to the complex product
structure, but they are always Ricci-flat, due to Corollary
\ref{ricci-flat}.

All these nilpotent Lie algebras have rational structure
constants, so that, by a theorem of Malcev, the corresponding
simply connected nilpotent Lie groups admit lattices, and hence
the associated nilmanifolds admit hypercomplex structures
invariant by the action of the group. If we take one of these
nilmanifolds with hypercomplex structure such that $\nhc$ is non
flat, then we obtain a hypercomplex nilmanifold with holonomy
contained in $SL(n,\HH)$, since the Ricci tensor vanishes. It was
proved in \cite{V} that the holonomy of the Obata connection of an
abelian hypercomplex structure on a nilmanifold is always
contained in $SL(n,\HH)$ where $4n$ is the dimension of the
nilmanifold. Here we can produce examples of hypercomplex
nilmanifolds with holonomy of the Obata connection also contained
in $SL(n,\HH)$ (since it is Ricci-flat) but with non abelian
complex structures. Note that none of these hypercomplex
structures underlies a hyperK\"ahler structure on the nilmanifold,
since there are no K\"ahler metrics on nilmanifolds, except for
the torus.

\smallskip

\begin{exam}
Consider the Lie algebra $\ggo=(0,0,0,12,13,14)$; it has a basis
$\{e_1,\ldots,e_6\}$ such that \[ [e_1,e_2]=-e_4,\quad
[e_1,e_3]=-e_5,\quad [e_1,e_4]=-e_6.\] According to Theorem
\ref{nonflat}, $\ggo$ admits both flat and non flat complex
product structures. An example of the former kind is given by
$\{J,\,E_1\}$, with
\[ Je_1=-e_2,\quad Je_3=-e_4,\quad Je_5=-e_6; \quad E_1|_{\ggo_+}=\Id,\quad
E_1|_{\ggo_-}=-\Id,\] where $\ggo_+$ and $\ggo_-$ are the Lie
subalgebras of $\ggo$ given by
\[ \ggo_+=\text{span}\{e_1,e_3,e_5\},\quad
\ggo_-=\text{span}\{e_2,e_4,e_6\},\] while an example of the
latter kind is given by $\{J,\,E_2\}$, with the same $J$ as above
and $E_2|_{\ggo'_+}=\Id,\quad E_2|_{\ggo'_-}=-\Id$, where
$\ggo'_+$ and $\ggo'_-$ are the Lie subalgebras of $\ggo$ given by
\[ \ggo'_+=\text{span}\{e_1,e_4,e_6\},\quad
\ggo'_-=\text{span}\{e_2,e_3,e_5\}.\] From Theorem 3.3 in
\cite{AS}, it follows that each complex product structure
$\{J,\,E_k\}$ gives rise to a hypercomplex structure
$\{\JJ,\II_k\},\, k=1,2$, on the Lie algebra $\hat{\ggo}=
(\ggo^\CC)_\RR$. This Lie algebra has a basis
$\{e_1,\ldots,e_6,\hat{e}_1,\ldots,\hat{e}_6\}$ such that
\begin{gather*}
[e_1,e_2]=-[\hat{e}_1,\hat{e}_2]=-e_4,\quad
[e_1,e_3]=-[\hat{e}_1,\hat{e}_3]=-e_5,\quad
[e_1,e_4]=-[\hat{e}_1,\hat{e}_4]=-e_6,\\
[\hat{e}_1,e_2]=[{e}_1,\hat{e}_2]=-\hat{e}_4,\quad
[\hat{e}_1,e_3]=[{e}_1,\hat{e}_3]=-\hat{e}_5,\quad
[\hat{e}_1,e_4]=[{e}_1,\hat{e}_4]=-\hat{e}_6.
\end{gather*}
The hypercomplex structure $\{\JJ,\II_1\}$ on $\hat{\ggo}$ is
given by
\begin{gather*}
\JJ e_1=-e_2,\quad \JJ e_3=-e_4,\quad \JJ e_5=-e_6,\quad
\JJ\hat{e}_1=-\hat{e}_2,\quad \JJ\hat{e}_3=-\hat{e}_4,\quad
\JJ\hat{e}_5=-\hat{e}_6,\\
\II_1e_1=\hat{e}_1,\quad \II_1e_3=\hat{e}_3,\quad
\II_1e_5=\hat{e}_5,\quad \II_1e_2=-\hat{e}_2,\quad
\II_1e_4=-\hat{e}_4,\quad \II_1e_6=-\hat{e}_6,
\end{gather*}
while, on the other hand, the hypercomplex structure
$\{\JJ,\II_2\}$ on $\hat{\ggo}$ is given by
\[
\II_2e_1=\hat{e}_1,\quad \II_2e_4=\hat{e}_4,\quad
\II_2e_6=\hat{e}_6,\quad \II_2e_2=-\hat{e}_2,\quad
\II_2e_3=-\hat{e}_3,\quad \II_2e_5=-\hat{e}_5,\] and $\JJ$ as
above. Due to the choice of the complex product structures, the
Obata connection associated to $\{\JJ,\II_1\}$ is flat, whereas
the Obata connection associated to $\{\JJ,\II_2\}$ is non-flat but
Ricci-flat. Note that these hypercomplex structures are not
abelian, since the Lie algebra $\ggo$ admits no abelian complex
product structure (see \S\ref{conc}).
\end{exam}

\

\section{Example: The 3-dimensional complex Heisenberg Lie group}
We consider the Lie algebra $\ggo=(0,0,0,0,13+42,14+23)$, which
has a basis $\{e_1,\ldots,e_6\}$ such that
$[e_1,e_3]=[e_,e_2]=-e_5,\, [e_1,e_4]=[e_2,e_3]=-e_6$. This Lie
algebra admits a matrix realization as $(3\times 3)$ complex
matrices in the following way:
\[ \ggo=\left\{\begin{pmatrix} 0&z_1&z_3\cr 0&0&z_2\cr 0&0&0\cr\end{pmatrix} : z_1,z_2,z_3\in\CC\right\},\]
where we can identify
\[ e_1=E_{1,2},\;e_2=ie_1,\;e_3=E_{2,3},\;e_4=ie_3,\;e_5=-E_{1,3},\;e_6=ie_5,\]
with $E_{r,s}$ the $(3\times 3)$ complex matrix whose only non
zero element is in the $(r,s)$ position and is equal to $1$. The
corresponding 6-dimensional simply connected nilpotent Lie group
$G$ can be described as a matrix group as follows:
\[ G=\left\{\begin{pmatrix} 1&z_1&z_3\cr 0&1&z_2\cr 0&0&1\cr\end{pmatrix} : z_1,z_2,z_3\in\CC\right\}.\]
That is, $G$ is the 3-dimensional complex Heisenberg Lie group,
considered as a real Lie group.

\smallskip

$\ri$ Let us consider the complex product structure $\{J,E\}$ on
$\ggo$ given by
\[ Je_1=e_3,\quad Je_2=e_4,\quad Je_5=e_6;\quad E|_{\ggo_+}=\Id,\quad E|_{\ggo_-}=-\Id,\]
where $\ggo_+$ and $\ggo_-$ are the Lie subalgebras of $\ggo$
given by
\[ \ggo_+=\text{span}\{e_1,e_2,e_5\},\quad \ggo_-=\text{span}\{e_3,e_4,e_6\}.\]
Note that both subalgebras are abelian, so that $\{J,E\}$ is an
abelian complex product structure on $\ggo$. The connected
subgroups of $G$ corresponding to $\ggo_+$ and $\ggo_-$ are given
by
\[ G_+=\left\{\begin{pmatrix} 1&z_1&x_3\cr 0&1&0\cr 0&0&1\cr\end{pmatrix} : z_1\in\CC,x_3\in\RR\right\},\quad
G_-=\left\{\begin{pmatrix} 1&0&iy_3\cr 0&1&z_2\cr
0&0&1\cr\end{pmatrix} : z_2\in\CC,y_3\in\RR\right\}.\] Since $G$
is nilpotent and simply connected, the subgroups $G_+$ and $G_-$
are also simply connected and closed in $G$. Moreover, it is easy
to verify in this case that $G=G_+\cdot G_-$, i.e., the
multiplication on $G$ defines a diffeomorphism $G_+\times G_-\lra
G,\,(g_+,g_-)\mapsto g_+g_-$; therefore, $(G,G_+,G_-)$ is a global
double Lie group, with $G_+\cong\RR^3\cong G_-$ (see \cite{LW}).

\smallskip

$\rii$ Let us consider now the complex product structure $\{J,E\}$
on $\ggo$ given by
\[ Je_1=-(e_4-e_2),\quad Je_2=-(e_1+e_3),\quad Je_5=e_6;\quad E|_{\ggo_+}=\Id,\quad E|_{\ggo_-}=-\Id,\]
where $\ggo_+$ and $\ggo_-$ are the Lie subalgebras of $\ggo$
given by
\[ \ggo_+=\text{span}\{e_1,e_2,e_5\},\quad \ggo_-=\text{span}\{e_1+e_3,e_4-e_2,e_6\}.\]
Note that $\ggo_+$ is abelian (and is the same subalgebra that
appeared in the case $\ri$), whereas $\ggo_-$ is isomorphic to
$\hg_3$. The connected subgroup of $G$ corresponding to $\ggo_+$
is the same subgroup $G_+$ appearing in the case $\ri$, while, on
the other hand, the subgroup $G_-$ corresponding to $\ggo_-$ is
\[ G_-=\left\{\begin{pmatrix} 1&\overline{z_2}&\frac{1}{2}|z_2|^2+iy_3\cr 0&1&z_2\cr 0&0&1\cr\end{pmatrix} :
z_2\in\CC,y_3\in\RR\right\}.\] In this case again we have again
that $G=G_+\cdot G_-$, so that $(G,G_+,G_-)$ is also a global
double Lie group, with $G_+\cong\RR^3,\, G_-\cong H_3$.

$\riii$ Let us consider now the complex product structure
$\{J,E\}$ on $\ggo$ given by
\[ Je_1=-(e_4-e_2),\quad Je_2=-(e_1+e_3),\quad Je_5=e_6;\quad E|_{\ggo_+}=\Id,\quad E|_{\ggo_-}=-\Id,\]
where $\ggo_+$ and $\ggo_-$ are the Lie subalgebras of $\ggo$ given by
\[ \ggo_+=\text{span}\{e_1+e_2+e_3,e_1-e_2+e_4,e_5-e_6\},\quad \ggo_-=\text{span}\{-e_1+e_2-e_3,e_1+e_2-e_4,e_5+e_6\}.\]
Note that both subalgebras are isomorphic to $\hg_3$. The
connected subgroups of $G$ corresponding to $\ggo_+$ and $\ggo_-$
are given by
\begin{gather*}
G_+=\left\{\begin{pmatrix}
1&(1+i)\overline{z}&\left(\ft|z|^2+t)+i(\ft|z|^2-t\right)\cr
0&1&z\cr 0&0&1\cr\end{pmatrix} :
z\in\CC,\,t\in\RR\right\},\\
G_-=\left\{\begin{pmatrix}
1&(1-i)\overline{z}&\left(\ft|z|^2+t)-i(\ft|z|^2-t\right)\cr
0&1&z\cr 0&0&1\cr\end{pmatrix} : z\in\CC,\,t\in\RR\right\}.
\end{gather*}
In this case we have again that $G=G_+\cdot G_-$, so that
$(G,G_+,G_-)$ is also a global double Lie group, with both $G_+$
and $G_-$ isomorphic to $H_3$.

\

\noindent{\textbf{Acknowledgements.} The author would like to
thank Isabel Dotti for useful discussions and support during the
preparation of the paper.}

\


\begin{thebibliography}{99}

\bibitem{A}
Andrada, A.: Complex product structures and affine foliations.
Ann.~Global Anal.~ Geom.~{\bf 27} (2005), 377--405.

\bibitem{ABDO}
Andrada, A.; Barberis, M.~L.; Dotti, I.; Ovando, G.: Product
structures on four dimensional solvable Lie algebras. Homology
Homotopy Appl.~{\bf 7} (2005), 9--37.

\bibitem{AD}
Andrada, A.; Dotti, I.: Double products and hypersymplectic
structures on $\RR^{4n}$. Commun. Math.~Phys.~{\bf 262} (2006),
1--16.

\bibitem{AS}
Andrada, A.; Salamon, S.: Complex product structures on Lie
algebras. Forum Math.~{\bf 17} (2005), 261--295.

\bibitem{Au}
Auslander, L.: Simply transitive groups of affine motions.
Am.~J.~Math.~{\bf 99} (1977), 809--826.

\bibitem{BD}
Barberis, M.L.; Dotti Miatello, I.: Hypercomplex structures on a
class of solvable Lie groups. Quart.~J.~Math.~Oxford {\bf 47}
(1996), 389--404.

\bibitem{BDM}
Barberis, M.L.; Dotti Miatello, I.; Miatello, R.: On certain
locally homogeneous Clifford manifolds. Ann.~Global Anal.~Geom.~{\bf 13} (1995), 289--301.

\bibitem{B}
Benoist, Y: Une nilvari\'et\'e non affine. J.~Differential Geom.~{\bf 41}
(1995), 21--52.

\bibitem{BV}
Blazi\'c, N.; Vukmirovi\'c, S.: Four-dimensional Lie algebras with
a para-hypercomplex structure. Preprint, available at
arXiv:math.DG/0310180

\bibitem{BM}
Bromberg, S.; Medina, A.: A homogeneous space-time model
with singularities. J.~Math.~Phys.~{\bf 41} (2000), 8190--8195.

\bibitem{CL}
Camacho, C.; Lins Neto, A.: Geometric Theory of Foliations.
Birkh\"auser, Boston, Massachusetts, 1985.

\bibitem{CFGU1}
Cordero, L.; Fernandez, M.; Gray, A.; Ugarte, L.: Nilpotent
complex structures on compact nilmanifolds. Rend.~Circolo Mat.~Palermo {\bf 49} suppl.~(1997), 83--100.

\bibitem{CFGU2}
Cordero, L.; Fernandez, M.; Gray, A.; Ugarte, L.: Compact
nilmanifolds with nilpotent complex structure: Dolbeault
cohomology. Trans.~Amer.~Math.~Soc.~{\bf 352} (2000), 5405--5433.

\bibitem{CFU}
Cordero, L.; Fernandez, M.; Ugarte, L.: Abelian complex structures
on 6-dimensional compact nilmanifolds. Comment.~Math.~Univ.~Carolinae {\bf 43} (2002), 215--229.

\bibitem{DDI}
De Cat, T.; Dekimpe, K.; Igodt, P.: Translations in simply
transitive affine actions of Heisenberg type Lie groups. Linear
Algebra Appl.~{\bf 359} (2003), 101--111.

\bibitem{DF}
Dotti, I.; Fino, A.: Abelian hypercomplex 8-dimensional
nilmanifolds. Ann.~Global Anal.~Geom.~{\bf 18} (2000), 47--59.

\bibitem{DF2}
Dotti, I. and Fino, A.: Hypercomplex nilpotent Lie groups. In:
Global Differential Geometry: The Mathematical Legacy of Alfred
Gray, 310--314, Contemp. Math. {\bf 241}, Amer. Math. Soc., 2001.

\bibitem{F}
Fried, D.: Distality, completeness and affine structures.
J.~Differ.~Geom.~{\bf 24} (1986), 265--273.

\bibitem{FG}
Fried, D.; Goldman, W.: Three dimensional affine crystallographic
groups. Advances in Math.~{\bf 47} (1983), 1--49.

\bibitem{K}
Kim, H.: Complete left-invariant affine structures on nilpotent
Lie groups. J.~Differential Geom.~{\bf 24} (1986), 373--394.

\bibitem{KK}
Kaneyuki, S.; Kozai, M.: Paracomplex structures and affine
symmetric spaces. Tokyo J.~Math. {\bf 8} (1985), 81--98.

\bibitem{L}
Libermann, P.: Sur les structures presque paracomplexes.
C.~R.~Acad.~Sci.~Paris {\bf 234} (1952), 2517--2519.

\bibitem{LW}
Lu, J.-H. and Weinstein, A.: Poisson Lie groups, dressing
transformations and Bruhat decompositions. J. Diff. Geom. {\bf 31}
(1990), 501--526.

\bibitem{M}
Magnin, L.: Sur les alg\`ebres de Lie nilpotentes de dimension
$\leq 7$. J.~Geom.~Phys.~{\bf 3} (1986), 119--144.

\bibitem{Ma}
Malcev, A.I.: On a class of homogeneous spaces. Reprinted in
Amer.~Math.~ Soc.~Trans.~ Ser.~1 {\bf 9} (1962), 276--307.

\bibitem{S}
Salamon, S.: Complex structures on nilpotent Lie algebras. J.~Pure
Appl.~Algebra {\bf 157} (2001), 311--333.

\bibitem{Se}
Segal, D.: The structure of complete left-symmetric algebras.
Math.~Ann.~{\bf 293} (1992), 569--578.

\bibitem{V}
Verbitsky, M.: Hypercomplex manifolds with trivial canonical
bundle and their holonomy. Preprint, 2004: arXiv:math.DG/0406537.
\end{thebibliography}
\end{document}